\documentclass[twoside,12pt, leqno]{article}
\usepackage{amsmath,amscd,amsthm,amssymb,amsxtra,latexsym,epsfig,epic,graphics}
\usepackage[matrix,arrow,curve]{xy}
\usepackage{graphicx}
\usepackage{diagrams}
\usepackage[english]{babel}
\usepackage[alphabetic,lite]{amsrefs} % for bibliography

\usepackage[OT2,OT1]{fontenc}
\newcommand\cyr{%
\renewcommand\rmdefault{wncyr}%
\renewcommand\sfdefault{wncyss}%
\renewcommand\encodingdefault{OT2}%
\normalfont
\selectfont}
\DeclareTextFontCommand{\textcyr}{\cyr}

\def\cprime{\char126}

\usepackage{hyperref}
\hypersetup{colorlinks=true,backref=true,citecolor=blue}

\def\antiddot{\mathinner{\mkern1mu\raise1pt\vbox{\kern7pt\hbox{.}}\mkern2mu
        \raise4pt\hbox{.}\mkern2mu\raise7pt\hbox{.}\mkern1mu}}

\newcommand{\PP}{{\mathbb P}}

\newcommand{\ZZ}{{\mathbb Z}}

\newcommand{\HH}{{\rm{H}}}

\newcommand{\coker}{{\rm{coker}\,}}

\newcommand{\s}{\mathcal}

\newcommand{\sB}{{\s B}}

\newcommand{\sF}{{\s F}}
\newcommand{\sG}{{\s G}}

\newcommand{\sL}{{\s L}}

\newcommand{\sO}{{\s O}}

\newcommand{\cO}{{\s O}}
\newcommand{\cF}{{\s F}}

\newcommand{\tensor}{\otimes}

\newcommand{\punkt}{\hspace{-.3ex}\raise.15ex\hbox to1ex{\Huge.}}

\def \fix#1 {{\hfill\break \bf (( #1 ))\hfill\break}}

\DeclareMathOperator{\Sym}{Sym}
\DeclareMathOperator{\reg}{reg}

\DeclareMathOperator{\Hom}{Hom}

\DeclareMathOperator{\depth}{depth}

\newcommand{\gm}{\mathfrak m}

\newtheorem{theorem}{Theorem}[section]
\newtheorem{lemma}[theorem]{Lemma}
\newtheorem{proposition}[theorem]{Proposition}
\newtheorem{corollary}[theorem]{Corollary}

\theoremstyle{definition}
\newtheorem{definition}[theorem]{Definition}
\newtheorem{notation}[theorem]{Notation}

\newtheorem{remark}[theorem]{Remark}

\newtheorem{example}[theorem]{Example}

\newtheorem{algorithm}[theorem]{Algorithm}

\def\bF{{\bf F}}
\def\bR{{\bf R}}
\def\bL{{\bf L}}
\def\bT{{\bf T}}
\def\bU{{\bf U}}
\def\bW{{\bf BW}}

\def\b1{{\bf 1}}
\def\b1{{1^t}}
\def\o{{\emptyset}}
\def\AA{{\mathbb A}}

\def\PP{{\mathbb P}}
\def\P{{\mathbb P}}

\def\HHH{{\mathbb H}}

\def\image{{\rm image}}
\def\cE{{\cal E}}
\def\cF{{\cal F}}

\def\cornerT#1{{T_{\Rsh \kern -1pt #1}}}
\def\Tate#1#2#3#4{{T_#1(#2,#3,#4)}}
\def\fix#1{{\bf ***} #1 {\bf ***}}

\DeclareMathOperator{\lincplx}{{\rm lincplx}}
\DeclareMathOperator{\cplx}{{\rm cplx}}
\DeclareMathOperator{\grmod}{{\rm grmod}}
\DeclareMathOperator{\RHom}{{\rm RHom}}
\DeclareMathOperator{\Tail}{{\rm Tail}}

\def\PPn{{\PP^{n_1} \times \cdots  \times \PP^{n_t}}}

\makeatletter
\def\Ddots{\mathinner{\mkern1mu\raise\p@
\vbox{\kern7\p@\hbox{.}}\mkern2mu
\raise4\p@\hbox{.}\mkern2mu\raise7\p@\hbox{.}\mkern1mu}}
\makeatother

\newtheorem{question}[theorem]{Question}

\usepackage{times}
\newdimen\x \x=12pt

\usepackage{color}

\date{}
\title{Tate Resolutions\\ for Products of Projective Spaces\\
\vskip .2cm \small{\it for Ngo Viet Trung on the occasion of his sixtieth birthday}}
\author{David Eisenbud, Daniel Erman, and Frank-Olaf Schreyer
\footnote{This paper reports on work started during the Commutative Algebra Program, 2012-13,
at the Mathematical Sciences Research Institute in Berkeley, and continued at the Mathematisches Forschungsinstitut Oberwolfach. We are grateful
to these institutes for providing  a beautiful and exciting environment for this work. The first and second authors are grateful to the 
National Science Foundation, and the second and third author are grateful to the Simons Foundation for partial support during this period. This material is based upon work supported by the National Science 
Foundation under Grant No. 0932078 000, while the authors were in 
residence at the Mathematical Science Research Institute in Berkeley, 
California in 2012--2013.}}

\begin{document}

\maketitle

\begin{abstract} We describe the Tate resolution of a coherent sheaf or complex of coherent sheaves on
a product of projective spaces. Such a resolution makes explicit all the cohomology of all twists of the sheaf, including, for example, the multigraded module of twisted global sections, and also the Beilinson monads of all twists.  Although the Tate resolution is highly infinite, any finite number of components can be computed efficiently, starting either from a Beilinson monad or from a multigraded module. 
\end{abstract}

\section*{Introduction}

A complex of coherent sheaves $\sF$ on projective space may be specified in finite terms by giving a complex of graded modules $\mathbb M$, or by giving a Beilinson monad $\sB$, that is, a finite complex  written in terms
of a strong exceptional sequence of the vector bundles $\wedge^{i}U$, where $U$ is the universal rank $n$ sub-bundle.
%,  such that 
%$$
%H_{i}(\sB) 
%= 
%\begin{cases} 
%\sF & \hbox{if } i=0\\ 0 &\hbox{otherwise.}
%\end{cases}
%$$
The complex $\mathbb M$ is a convenient way of simultaneously specifying the complexes of global sections $H^{0}(\sF(d))$ of twists $\sF(d)$ for all sufficiently large $d$, while the Beilinson monad is a convenient way of specifying the hypercohomology of $\sF$ itself; but the Beilinson monads of twists $\sF(d)$ generally look quite different.
The \emph{Tate resolution} of $\sF$ is a way of packaging all the cohomology spaces and  the Beilinson monads of all twists of $\sF$ simultaneously. It is a doubly infinite exact complex of finitely generated free modules over the exterior algebra $E$ that is the Koszul dual of the homogeneous coordinate ring of projective space. Any finite number of terms can be computed  efficiently, in terms of $E$-free resolutions of finitely generated $E$-modules, from either a module of twisted global sections or a Beilinson monad for $\sF$. Tate resolutions in this case were treated in \cite{EFS} and \cite{ESW}. In many cases they yield the fastest algorithms for computing cohomology.

In this paper we will provide an analogous, efficiently computable, construction of a Tate resolution for finite complexes of coherent sheaves on products of projective spaces. A new feature, which makes this case much more difficult than the case of a single projective space, is that there are no finitely generated modules among the terms of the Tate resolution. Despite this, we can
use the Tate resolution to effectively compute the monads and (hyper)cohomology of any finite number of twists of $\sF$ in terms of free resolutions of certain finitely generated modules over an appropriate exterior algebra.

To state the main results we introduce some notation:

Let $\PP=\PPn =  \PP(W_{1})\times\cdots\times \PP(W_{t})$ be a product of $t$ projective spaces over an arbitrary field $K$. Set $V_{i} = W_{i}^{*}$ and  $V = \oplus_{i}V_{i}$. Let $E$ be the $\ZZ^{t}$-graded exterior algebra on $V$, where elements of $V_{i}\subset E$ have degree $(0,\dots,0, -1, 0,\dots,0)$ with $-1$ in the $i$-th place.

Let $\omega_{E}$ be the free $E$-module of rank 1 with generator in multidegree
$(n_{1}+1, \dots, n_{t}+1)$; this module has socle in multidegree $0$ and is the injective
hull of the residue field $K$. We will generally write free complexes of 
$E$-modules as sums of twists of $\omega_{E}$. We have $\omega_E
=\Hom_K(E,K)$, where $K$ is regarded as a 1-dimensional vector space
concentrated in degree 0.

A $\ZZ^{t}$-graded complex of $E$-modules is called \emph{locally finite} if 
the sum of the graded components of each multidegree is finite dimensional. 

Let $U_k = \ker(H^0(\PP^{n_k},\sO(1))\tensor \sO \to \sO(1))$ denote  the tautological subbundle on  $\PP^{n_k}$ of rank $n_k$. For $a\in \ZZ^{t}$ we set
$$
U^a := \boxtimes_{k=1}^t \Lambda^{a_k} U_k=\pi_1^*\Lambda^{a_1}U_1 \otimes \cdots \otimes \pi_t^*\Lambda^{a_t} U_t
$$
the tensor product of the pullbacks to $\PP$ of exterior powers of the $U_k$. 
Of course $U^a$ is nonzero if only if $0 \le a \le n=(n_1,\ldots,n_t)$, where the partial order on multi-indices is termwise. The $U^{a}$ form a \emph{full strong exceptional collection} for the derived category $D^{b}(\P)$~\cite[Def. 8.31]{huybrechts}. In particular every sheaf can be expressed as the homology of a complex whose terms are direct sums of the sheaves $U^{a}$, called a \emph{Beilinson monad} for $\sF$.
%(~\cite[Lemma~1.6]{kapranov}).

Consider the additive functor on the category of direct sums of finitely generated free graded right $E$-modules defined on objects by
$$
\bU\colon \omega_E(a) \mapsto U^a.
$$
We define $\bU$ on morphisms by using the identifications
$$ \Hom_{E}(\omega(a), \omega(b)) = E_{b-a} = 
\Hom_{\P}(U^{a}, U^{b}).
$$
This identifies the morphism induced by multiplication by an element 
$$
e\in E_{b-a}= \otimes_{k=1}^t \Lambda^{a_k-b_k}V_k
$$ 
with a morphism of sheaves
$$
\bU\left(\omega_{E}(a) \rTo^{ e\wedge-} \omega_{E}(b)\right) =
\left(U^{a}\rTo^{e\neg } U^{b}\right),
$$
where 
$e\neg $ is the map 
%from 
%$$
%U^{a}\subset \otimes_{k=1}^t \Lambda^{a_k}W_k\otimes_K \sO
%$$
%to 
%$$
%U^{b}\subset \otimes_{k=1}^t \Lambda^{b_k}W_k \otimes_K \sO.
%$$
induced by the contraction operator, which we write in the same way,
$$
\otimes_{k=1}^t \Lambda^{a_k}W_k\otimes_K \sO\rTo^{e\neg}
\otimes_{k=1}^t \Lambda^{b_k}W_k \otimes_K \sO.
$$
See for example \cite{EFS} for the case $t=1$.

%as indicated in the diagram
%$$
%\begin{diagram}
% \otimes_{k=1}^t \Lambda^{a_k}W_k\otimes_K \sO&\rTo^{\neg e}&
%\otimes_{k=1}^t \Lambda^{b_k}W_k \otimes_K \sO\\
%\uTo&&\uTo\\
%U^{a}&\rTo^{\neg e}& U(b)
%\end{diagram}
%$$
%$$
%\xymatrix{
%\omega_E(a)  \ar[d]_{\times e}& \mapsto & U^a \ar[d]^{\neg e} & \hookrightarrow & \otimes_{k=1}^t \Lambda^{a_k}W_k\otimes_K \sO \ar[d] ^{\neg e}\cr
%\omega_E(b) &\mapsto & U^b & \hookrightarrow & \otimes_{k=1}^t \Lambda^{b_k}W_k \otimes_K \sO \cr }. 
%$$
%Applying $\bU$ to a locally finite complex of free $E$-modules,
% we obtain a bounded complex
%$
%\bU(\T).
%$

If $T$ is a multigraded complex of free $E$-modules,  $I\subsetneq \{1,\dots,t\}$
is a proper subset,
and $c = (c_{i})_{i\in I}$ are integers,
then the \emph{$I$-th strand of $T$} through $c$
is the subquotient complex obtained from $T$ by taking all the
free summands of terms of $T$ of the form $\omega_E(a)$ where $a_{i} = c_{i}$ for all 
$i\in I$. When all the $c_{i}$ are zero, we speak simply of the $I$-th strand of $T$.
Thus for example $T$ itself is the $\emptyset$-strand. 

\goodbreak

We will say that a locally
finite $\ZZ^{t}$-graded complex of free $E$-modules is a 
 \emph{Tate Resolution} if, for every
multi-index $c$, all strands of $T$ through $c$ are exact. A complex
$T$ of free $E$-modules is called 
minimal if $T\otimes_{E}K$ has zero differential. Every complex
of free $E$-modules
is isomorphic to the direct sum of a minimal complex and a split exact complex.

\begin{theorem}\label{Tate existence}
  For any sheaf $\sF$ on $\PP$ there is a unique minimal Tate
  resolution $\bT(\sF)$ such that 
$\bU\bigl(\bT(\sF)(c)[|c|]\bigr)$ is a Beilinson monad for $\sF(c)$ for every $c\in \ZZ^{t}$.
Moreover, 
$$
\bT(\sF)^d = \oplus_{a\in \ZZ^t} \Hom_K(E,H^{d-|a|}(\PP, \sF(a))),
$$
where the cohomology $H^{d-|a|}(\PP,\sF(a))$ is regarded as a vector
space concentrated in degree $d-|a|$.
\end{theorem}

The second statement follows from the first using the well-known
result of Theorem~\ref{Terms of the Beilinson monad}. The construction
of the Tate resolution is given in Section~\ref{construction}, and the
proofs of its properties are given in Sections~\ref{monads} and \ref{exactness properties}.

\begin{example} Consider the case where $t=1$ and $\PP = \PP^1$.
The Tate resolution of $\sO_{\PP^{1}}$ on $\PP^{1}$ is a complex of the form
$$
 T:=\bT(\sO_{\PP^{1}}) = \cdots \rTo \omega_{E}^{2}(3) \rTo^{
\begin{pmatrix}
e_{0}&e_{1} 
\end{pmatrix}
} \omega_{E}(2)\rTo^{
\begin{pmatrix}
 e_{0}e_{1}
\end{pmatrix}}
 \omega_{E}\rTo^{
\begin{pmatrix}
e_{0}\\ e_{1} 
\end{pmatrix}
} \omega_{E}^{2}(-1) \rTo \cdots.
$$
Every term in the resolution $\bT(\sF)$ is finitely generated in this case; indeed,
this is true whenever $t=1$.

By contrast, if we take $t>1$ then each term in the Tate resolution of
\emph{any} non-zero sheaf is infinite. For instance, in the case
$t=2$ with $\PP = \PP^1\times \PP^1$
and $\sF = \sO_{\PP^{1}\times\PP^{1}} = \sO_{\PP^{1}}\boxtimes
\sO_{\PP^{1}}$,
then $\bT(\sF) = T\boxtimes T$, a complex in which \emph{every} term
is infinite. By the formula in Theorem~\ref{Tate existence} we have, for example,
\begin{align*}
\bT(\sF)^0 &= 
\oplus_{(p,q)\in \ZZ^2} \Hom_K\left(E,H^{0-p-q}\left(\PP, \sF(p,q)\right)\right)\\
&= \Hom_K\left((E,H^{0}(\PP, \sF)\right) \oplus \bigoplus_{p\in \ZZ} \Hom_K\left((E,H^1(\PP, \sF(p,-p-1))\right)\\
&= H^{0}(\PP, \sF)\otimes \omega_E \oplus \bigoplus_{p\in \ZZ} H^1\left(\PP, \sF(p,-p-1)\right)\otimes \omega_{E}(-p,p+1) \\
&= \omega_{E} \oplus\bigoplus_{p\in \ZZ}\omega_{E}(-p,p+1)^{p(p+1)}
\end{align*}
Other examples are given in Section~\ref{P1xP1 examples}.
\end{example}

\bigskip

Since $\bT(\sF)$ is locally finite we can form, for any finite interval 
$A$ in $\ZZ^{t}$, the finite subquotient complex $T'$ of $\bT(\sF)$ with terms 
$$
T'^d \cong  \bigoplus_{ a\in A} 
\Hom_K(E, H^{d-|a|}(\sF(a)).
$$
We give an algorithm, implemented in Macaulay2, for computing any such finite subquotient complex~\cite{TateOnProducts}.
The algorithm makes use of free resolutions over $E$, and can be executed starting either from a multigraded complex of modules representing
the global sections of high twists of $\sF$ or from a Beilinson monad for $\sF$.
In particular, our method computes any finite number of cohomology groups of a coherent sheaf on $\PPn$ without passing to a Segre embedding, a process that would introduce a much larger number of variables.

We can use the Tate resolution of $\sF$ to compute the direct image complex  along the projection
$\pi_{J}: \P\to \P_{J}:= \prod_{j\in J}\P^{n_{j}}$ for any  subset $J\subset \{1,\dots,t\}$. 
Let $I$ be the complement of $J$. The differential on the $I$-th strand of $\bT(\sF)$ is given by matrices with entries in $E_{J}:=\wedge \oplus_{j\in J}V_{j}$.
Thus the $I$-th strand of $\bT(\sF)$  has the form $T_{J}\otimes \omega_{E_{I}}$, where $T_{J}$ is a complex over $E_{J}$.
Because a strand of a strand is a strand, $T_{J}$ is  a Tate resolution on $\PP_{J}$.

\begin{corollary}\label{strands}
With notation as above, $\bU(T_{J})$ is a Beilinson monad for $R\pi_{J*}\sF$.
\end{corollary}

Since $\bU(\bT(\sF))$ is quasi-isomorphic to $\sF$, we see that
$\bU(\bT(\sF(c)))$ must be quasi-isomorphic to $\bU(\sF)(c)$. The
formula for $\bT(\sF(c))$ above suggests the following:

\begin{theorem}[Twist and shift formula]\label{twist and shift}
Let $T$ be a Tate resolution. Then
$
\bU(T(c))
$
and
$\bU(T)(c)[-c]$
are quasi-isomorphic. In particular, for any non-empty subset $J \subset \{1,\ldots t \}$ with complement $I$, the $I$-th strand of $T$ through $c$ computes the direct image
$R{\pi_J}_*\bigr( \bU(T)(c)\bigl)$ along $\pi_J: \PP \to \PP_J$.
\end{theorem}

\medskip
\noindent
\begin{example}
The quasi-isomorphism $\bU(\bT(\sF(c)))\simeq\bU(\sF)(c)$
above is generally \emph{not} an isomorphism. For example
 $\bU(\bT(\sO_{\PP^{1}\times\PP^{1}}))$ is the complex 
$$
0 \to \sO_{\PP^{1}\times\PP^{1}} \to 0 .
$$
By contrast, $\bU(\bT(\sO_{\PP^{1}\times\PP^{1}}(3,-2)))$
is 
$$
0 \to (U_1 \boxtimes U_2)^{\oplus 6} \to (U_1 \boxtimes \sO)^{\oplus 3} \oplus (\sO \boxtimes U_2)^{\oplus 8} \to (\sO \boxtimes \sO)^{\oplus 4} \to 0
$$
which we may rewrite as
$$
0 \to \sO^6(-1,-1) \to \sO^3(-1,0) \oplus \sO^8(0,-1) \to \sO^4 \to 0.
$$
Thus $\bU(\bT(\sO_{\PP^{1}\times\PP^{1}}))(3,-2)$ and $\bU\bigr(\bT(\sO_{\PP^{1}\times\PP^{1}})(3,-2)[-1]\bigl)$ are  only quasi-isomorphic.
\end{example}
\bigskip

Our methods can also be used to treat the hypercohomology of complexes of
sheaves. We say that a complex $U$ is a \emph{Beilinson
  representative} for a bounded complex $\sF$ of sheaves if $\sF$ is
quasi-isomorphic to $U$ and each term of $U$ is a direct sum of copies
of the sheaves $U^a$. Our construction of $\bT(\sF)$ when $\sF$ is a
sheaf generalizes immediately to the case when $\sF$ is a bounded
complex of sheaves, and then $\bU(\bT(\sF))$ is a Beilinson
representative of $\sF$. As above, it follows that 
\begin{equation}\label{eqn:hypercohomology}
\bT(\sF)^{d} \cong \bigoplus_{a \in \ZZ^t }  \Hom_K\left(E, \HHH^{d-|a|}(\PP,\sF(a))\right).
 \end{equation}
where $\HHH^{d-|a|}$ denotes hypercohomology. For simplicity we will generally focus on the case
of sheaves, leaving some details of the case of complexes to the interested reader.

The construction of the Tate resolution is given in Section~\ref{construction}, and the connection with Beilinson monads as well as the proof of Theorem~\ref{Tate existence} appears in Section~\ref{monads}.

Associated to a Tate resolution are many exact  complexes built from induced subquotient complexes---see Theorem~\ref{exactness of T}; one type, the ``corner complexes'' defined in Section~\ref{exactness properties}, are the key to our finitistic construction of Tate resolutions and also to our proof of 
Corollary~\ref{strands} and the general twist and shift formula stated in Theorem~\ref{twist and shift}.

Given a finite complex $BW$ of free $E$-modules we can apply the functor
$\bU$ to get a finite complex of sheaves, and then apply the functor
$\bT$ to get a doubly infinite exact complex $T$ over $E$. In
Section~\ref{Tate from Monad} we show how to go directly from $BW$ to $E$ by a
computation in terms of free and injective resolutions over $E$.

In Section~\ref{splittings} we consider a broad generalization, to products of
projective spaces, of
Horrocks criterion for the splitting of a vector bundle. We prove that the criterion holds
under an additional hypothesis, which may not be necessary.

We are grateful to Mike Stillman and Dan Grayson, the authors of the
computer algebra system Macaulay2~\cite{M2}, without which we would
not have discovered the results in this paper! We are also grateful
to Christine Berkesch and Florian Geiss for useful conversations.

\section{Construction of the Tate resolution}\label{construction}
As above, $\PP=\PPn =  \PP(W_{1})\times\cdots\times \PP(W_{t})$ is a product of $t$ projective spaces over an arbitrary field $K$.  We let $V_{i} = W_{i}^{*}$, $W=\oplus_i W_i$ and  $V = \oplus_{i}V_{i}$. We let $S=\Sym(W)$ be the Cox ring of $\PP$, with the $\ZZ^t$ grading where elements in $W_i\subset S$ have degree $(0,\dots,0,1,0,\dots,0)$ with the $1$ in the $i$th place.  We also let $E$ be the exterior algebra on $V$ with the dual grading, i.e. elements of $V_{i}\subset E$ have degree $(0,\dots,0, -1, 0,\dots,0)$.  We write $\{x_{i,j}\}$ and $\{e_{i,j}\}$ for dual bases of $W_{i}$ and $V_{i}$.

By~\cite{cox}, there is a correspondence between finitely generated, multigraded $S$-modules and coherent sheaves on $\PP^n$.  More precisely, if $M$ is any finitely generated, multigraded $S$-module, then there is a corresponding coherent sheaf $\widetilde{M}$ on $\PP$, and every coherent sheaf arises in this way.  
%Moreover, $M$ and $M'$ yield the same sheaf if they agree up to $B$-torsion, where $B:=\cap_{i=1}^t (x_{i,0},\dots,x_{i,n_i})\subset S$ is the irrelevant ideal of the Cox ring of $\PP$. Thus, for instance, if $b\in \ZZ^t$ is some multidegree then $\widetilde{M_{\geq b}}=\widetilde{M}$.

%Let $\omega_{E}$ be the free $E$-module of rank 1 with generator in multidegree
%$(n_{1}+1, \dots, n_{t}+1)$; this module has socle in multidegree $0$ and is the injective
%hull of the residue field $K$. We will generally write free complexes of 
%$E$-modules as sums of twists of $\omega_{E}$.

%A $\ZZ^{t}$-graded complex of $E$-modules is called \emph{locally finite} if 
%the sum of the graded components of each multidegree is finite dimensional. 
%We will use the notations $S, E, \omega_{E},\dots$ of the introduction throughout the paper. We write
%$\{x_{i,j}\}$ and $\{e_{i,j}\}$ for dual bases of $W_{i}$ and $V_{i}$

%Throughout the paper $S$ and $E$ will denote the $\ZZ^{t}$-graded polynomial ring and exterior algebra
%$$
%S := \Sym \left(\oplus_{i=1}^{t} W_{i}\right); \quad E:= \wedge\left(\oplus_{i=1}^{t} V_{i}\right),
%$$
%where the $W_{i}$ and $V_{i}$ are finite-dimensional $K$-vector spaces, and
% $V_{i}$ is the vector space dual of $W_{i}$. 
 
 If $N$ is a $K$-vector space then we regard
 $\Hom_K(E,N)$  as a right $E$-module by the formula $(\phi e)(f) = \phi(ef)$.

%We first sketch the construction of the Tate resolution of a coherent sheaf, and then fill in the details.
We will use the BGG correspondence, which we briefly recall. 
\begin{definition}\label{def of R and L}
Let $M=\sum_{a\in \ZZ^t} M_a$ be a multigraded $S$-module. For $d\in \ZZ$ set $M_d = \sum_{a,|a|=d} M_a$, the sum of the components in total degree $d$. Let $\bR(M)$ be the complex
$$
\bR(M): \ldots \to \Hom_K(E,M_d) \to \Hom_K(E,M_{d+1}) \to \ldots 
$$
where $\Hom_K(E,M_d) = \omega_{E}\otimes_{K}M_{d}=\bR(M)^d$ is in cohomological degree $d$ and the differential is given by 
$$ 
\phi \mapsto \{ f \mapsto \sum x_{i,j} \phi(fe_{i,j}) \}.
$$ 
This is the same formula as in the singly graded case, but now $\bR(M)$ is a multigraded complex of free right $E$-modules.
Similarly, if $P=\sum_{a \in \ZZ^t} P_a$ is a multigraded right $E$-module, we define a complex of free $S$-modules, 
$$\bL(P): \ldots \to S\tensor_K P_d \to  S \tensor_K P_{d-1} \to \ldots$$
with $S\tensor_K P_d=\bL(P)_d$  in homological degree $d$ and $P_d= \sum_{a,|a|=d} P_a$. The differential is given by 
%\frank{we want R(M) to be a sequence of right E-modules, so that the matrices in Macaulay2 make sense, check once more left right.}
$$ s \tensor p \mapsto \sum sx_{i,j} \tensor  pe_{i,j}.$$

\end{definition} Note that our convention is $\deg e_{i,j} = -1 = -\deg x_{i,j}$.

The functor  $\bL$ defines an equivalence between the category $\grmod(E)$ of finitely generated  multigraded $E$-modules  and the category $\lincplx(S)$ of finite linear  complexes over $S$.
Indeed, suppose that
$$
\bF: \ldots \rTo^{\partial} F_d \rTo^{\partial} F_{d-1} \rTo^{\partial} \ldots
$$
is a linear complex with $F_d =\sum_{a: |a|=d} B_a \tensor_K S(-a)$, where we think of $B_a$ as a vector space in degree $0$. Let
$P = \sum_{d\in \ZZ} \sum_{a: |a|=d} P_a$ with $P_a=B_a(-a)$. We give $P$ the structure of a graded right $E$-module by setting
$ pe_{i,j} = p_{i,j}$  if
%
%act as
%$$P_a \ni p \mapsto p_{i,j} \in P_{a-\deg x_{i,j}}$$
%where 
$\partial (p\tensor 1) = \sum p_{i,j}\tensor x_{i,j}$.
 Because $\partial^2=0$,  this action extends to 
an action of the exterior algebra $E$, and $F = \bL(P)$.  
Similarly, $\bR: \grmod(S) \to \lincplx(E)$ is an equivalence between the category of multigraded $S$-module and  the category of linear free $E$-complexes.
Sometimes it is more convenient to index $\bL(P)$ cohomologically. For that purpose we define  $P^d=P_{-d}$.

 The functors $\bR$ and $\bL$ extend naturally to functors
 $$
\xymatrix{ \cplx(S) \ar[r]^{\bR} & \cplx(E) \ar[l]^{\bL} \cr}
$$
between the categories of complexes. We first define
$\bL(P[k]) = \bL(P)[k]$ and  $\bR(M[k])=\bR(M)[k]$, 
where, for any complex $F$ and $k\in \ZZ$ we define the shifted complex $F[k]$ by $F[k]^d= F^{k+d}$. Thus the shift operator commutes with the functors $\bR$ and $\bL$.  With this
convention we define $\bR$ on a complex of $S$-modules  by applying $\bR$ to each term,
and  taking the total complex of the resulting double complex; we define $\bL$ on a complex of $E$-modules similarly.

The shift operators also commute with the twist operators: $M[k](a) = M(a)[k]$ for any graded module over $S$ or $P$. However, the twist operator only commutes with $\bR$ and $\bL$ up to a shift, as in the following Lemma.
We can only shift a complex by an integer, not a multi-index so,  to simplify notation, we will usually write $M[b]$ in place of $M[|b|]$ (recall that $|b|$ denotes $\sum_{i}b_{i}$, \emph{not} the absolute value).

\begin{lemma}\label{TwistShift}
For any multigraded $S$-module  $M$, any multigraded $E$-module $P$,
and any $b\in \ZZ^{t}$, we have:
\begin{align*}
 \bR(M(b)) &= \bR(M)(b)[-b]\\
\bL(P(b)) &= \bL(P)(b)[-b].
\end{align*}
\end{lemma}
%$$
%(\bR(M))(b) = \bR(M[|b|](b))
%$$
%\david{I'd prefer to write this as
%$$
%\bR(M(b)) = \bR(M)(b)[|-b|]
%$$
%and similarly
%$$
%\bL(P(b)) = \bL(P)(b)[-b].
%$$
%}
%for any multigraded $S$-module  $M$ and any $b\in \ZZ^{t}$. 
%$$
%\bR(M(b)) = \bR(M)(b)[-b]
%$$
\begin{proof}
We regard  $M_a$ as a vector space in degree $a$, so that $M(b)_a=M_{b+a}(b)$. Thus
\begin{align*} 
\bR(M(b))^d  &= \sum_{a:|a|=d} \Hom_K(E,M_{b+a}(b)) \\
&= \sum_{c:|c|=d+|b|} \Hom_K(E,M_{c})(b) \\ &=(\bR(M)(b)[-b])^d 
\end{align*} 
The computation for $\bL$ is similar.
\end{proof}
%Similarly,
%%$$
%\begin{equation}\label{TwistShift}
%\bL(P[b](b)) = (\bL(P))(b).
%\end{equation}
%%$$
\medskip

A fundamental reciprocity result proven for the case $t=1$ in \cite{EFS}, Theorem 3.7, also holds in general:

\begin{theorem}[Reciprocity]\label{reciprocity} Let $M$ be a finitely generated graded $S$-module, and $P$ a finitely generated graded $E$-module.
Then $\bR(M)$ is an injective resolution of $P$ if and only if $\bL(P)$ is a projective resolution of $M$.
\end{theorem}

\begin{proof} We may regard $M$ and $P$ also as singly graded modules for which the Reciprocity Theorem~\cite[Theorem~3.7]{EFS}
applies. The constructions respect the finer grading.
\end{proof}

\begin{corollary}\label{quadrant Tate}
Let $M$ be a finitely generated multigraded $S$-module. $M_{\ge c}(c)$ has a linear resolution $F$ if and only if  $\bR(M_{\ge c})$ is acyclic. Moreover,  if we write $F=\bL(P)$ ,  then  $\HH^{|c|}(\bR(M_{\ge c})) \cong P(-c)$. 
\end{corollary}

\begin{proof} 
$\HH^0(\bR(M_{\ge c}(c))) \cong P$ is equivalent to $\HH^{|c|}(\bR(M_{\ge c})) \cong P(-c)$.
\end{proof}

We next recall  the construction from~\cite{EFS} of the Tate resolution of a sheaf $\sF$ on $\PP^n$. Suppose that
$\sF$ is the sheafification of the $\ZZ$-graded $S$-module $M$.
\begin{enumerate}
\item Choose an integer $c$ sufficiently large that the module $M_{\geq c}(c)$ has a linear resolution
(this is equivalent to $c\geq \reg M$), and so that $M$ contains no submodule of finite length (this is
satisfied if $c> \reg M$). It follows that
 $$
 \bR(M_{ \ge c}): 0 \to T^c \to T^{c+1} \to \ldots 
 $$
  is acyclic by Theorem~\ref{reciprocity} and is minimal because no element of $M$ is annihilated by every linear form of $S$.
  
\item  Attach to $\bR(M_{\ge c})$  a minimal free resolution
$ \ldots \to T^{c-2} \to T^{c-1} \to \HH^c(\bR(M_{\ge c})) \to 0$
to obtain a doubly infinite complex 
$$
\bT(\sF) :  \ldots \to T^{d-1} \to T^d \to T^{d+1} \to \ldots.
$$
\end{enumerate}
We have
\begin{theorem}[\cite{EFS} Thm 4.1]
  The complex $\bT(\sF)$ depends only on the sheaf $\sF$. Moreover,
$$ 
T^d= \sum_{i=0}^n \Hom_K(E, \HH^i(\sF(d-i))).
$$
\end{theorem}

Steps 1 and 2 work also in multigraded setting:
Given a $\ZZ^t$-graded $S$-module $M = \oplus_{a \in \ZZ^t} M_a$  and a multidegree $c \in \ZZ^t$ we denote by
$$
M_{\ge c} = \oplus_{a  \ge c} M_a
$$
its truncation.  Here $a \ge c$ stands for the componentwise partial order on multidegrees.  We note that, for any $c$, the sheaves $\widetilde{M}$ and $\widetilde{M_{\geq c}}$ are isomorphic~\cite[Proposition~5.3.10]{CLS}.

We will show that, for a sufficiently large multidegree $c$, the complex $\bR(M_{\ge c})$ is acyclic and minimal.
Thus we may attach a free resolution of  $\HH^{|c|}(\bR(M_{\ge c}))$ to obtain a doubly infinite complex of (finitely generated) free $E$-modules. 
However it is no longer true  that the resulting complex encodes all of the cohomology groups of
twists of $\sF= \widetilde M$. Instead, as we shall see, it encodes only the cohomology groups
of the twists that are $\geq c$ and those that are $\leq c-(1,\dots, 1)$; these are only $2$ of the $2^t$ orthants of $\ZZ^t$.  

To get a complex that encodes all cohomology groups we must pass to a restricted inverse limit. Fix a large multidegree $b$ and consider all $c \ge b$. Let $\Tail_c(M)$ denote the projective resolution of $\HH^{|c|}(\bR(M_{\ge c}))$. It is easy to see that the complexes 
$\bR( M_{\ge c})$ form a directed system of subcomplexes $\bR(M_{\ge c}) \hookrightarrow \bR(M_{\ge b})$.
We will show that the tails form an inverse system of complexes via a sequence of epimorphisms
$\Tail_c(M) \twoheadrightarrow \Tail_b(M)$ defined when $c>b\gg 0$.

\begin{definition}
The Tate resolution of  $\sF = \widetilde M$ on the product $\PP$ of $t$ projective spaces is the restricted inverse limit
$$
\bT(\sF)= \left({\lim_{\infty\leftarrow c}}' \Tail_{c}(M)\right)[-t+1].
 $$
where ``restricted'' means that, thinking of
$\lim_{\infty\leftarrow c} \Tail_{c}(M)$ as a subcomplex of
$\prod_{c} \Tail_{c}(M)$, we take only those sequences of elements 
of bounded degree.
\end{definition}

To make sense of this definition we must define the maps $\Tail_c(M) \twoheadrightarrow \Tail_b(M)$. The first
step is to prove some properties of high truncations of multigraded modules that are standard for singly graded modules. Recall that for any $J \subset \{1,\ldots t\}$, we denote the projection 
$\P\to \P_{J}:= \prod_{j\in J}\P^{n_{j}}$
by $\pi_J$.

\begin{proposition}\label{high truncation}
Let $M$ be a finitely generated multigraded $S$-module and let $\sF=\widetilde M$ be the coherent sheaf on $\PP=\PPn$ represented by $M$. There exists a multidegree $b$ such that, for any multidegree $c\ge b$:

\begin{enumerate}
\item The truncated twisted module $M_{\ge c}(c)$ is generated in degree $0$ and has a linear resolution, i.e. the free modules $F_k = \oplus_a S^{\beta_{k,a}}(-a)$ in the minimal free resolution
$$(0 \leftarrow M_{\ge c}(c) \leftarrow) F_0 \leftarrow F_1 \leftarrow \cdots$$
of $M_{\ge c}(c)$ satisfies $\beta_{k,a} \not=0$ only if $k=|a|$ and $a\geq 0$, 

\item $ M_{c} = H^0(\PP, \sF(c))$ and $H^p(\PP, \sF(c))=0$ for $p>0$.  More generally,
for any $J \subset \{1,\ldots t\}$ and $p>0$,  we have  $R^p {\pi_J}_* \sF(c)=0$. 

\end{enumerate}
\end{proposition}

\begin{definition}\label{defn suff positive}
  We call a multidegree $b$  \emph{sufficiently positive for $M$} if it satisfies
the conditions of  Proposition~\ref{high truncation}.
\end{definition}

\begin{proof} 
Consider first the case $M=S$. Write $\sO = \sO_{\P}$ for the structure sheaf. We claim that $0$ is sufficiently positive for $S$. 

Denote by $\gm_i=\langle W_i \rangle \subset S$ the minimal primes of the irrelevant ideal in the Cox ring. For every $c \ge 0$ the ideal
$$
S_{\ge c}= \bigcap_{i}\gm_i^{c_i} = \prod_{i} \gm_i^{c_i}
$$ 
has a resolution which is the tensor product of the resolutions of the ideals
$\gm_i^{c_i}$. Since (up to twist) these resolution are linear, the tensor product resolution of $S_{\ge c}(c)$ is linear, i.e. the first assertion holds for $S$ and every $c \ge 0.$ Furthermore, 
we have $S= \sum_{a } H^0(\sO(a_1,\ldots,a_t))$ and $\sum_{a \ge -n} H^p(\sO(a_1,\ldots,a_t)) =0$ for $p \ge1$ by the K\"unneth formula.
Finally  for a pair of complementary index sets $J \cup I = \{ 1, \ldots, t \}$ and multidegree $c=(c_J,c_I) \ge 0$, the truncated  module of global sections $(\pi_J)_*\sO(c)=H^0(\PP_I,\sO(c_{I})) \tensor \sO(c_J)$ has a linear resolution since $c_{J} \ge 0$.

We now consider an arbitrary finitely generated multigraded $S$-module $M$. By the Hilbert syzygy theorem, $M$ has a finite free multi-homogeneous
resolution
$$
0 \leftarrow M \leftarrow G_0 \leftarrow G_1 \leftarrow \cdots \leftarrow G_N \leftarrow 0
$$
with  $N \le \sum (n_i+1)$. Write $G_k = \oplus S(-a)^{\beta_{k,a}}$ and  
set $b_i = \max \{ a_i \mid \exists \beta_{k,a} \not=0 \}$. 

We claim that $b=(b_1,\ldots b_t)$ is sufficiently positive for $M$. Let $c \ge b$. Each $(G_k)_{\ge c}(c)$ has a linear resolution because, for every summand $S(-a)$, the truncation
$$
S(-a)_{\ge c}(c)=  \gm_1^{c_1-a_1} \cap \ldots \cap \gm_t^{c_t-a_t}(c-a)
$$
has a linear resolution. An iterated mapping cone (see Section~\ref{resolution matters}) over the complex $(G_\bullet)_{\ge c}(c)$ yields
a non-minimal resolution $F'$ with graded Betti numbers $\beta_{k,a}(F')$ non-zero only if $|a| \le k$.
Since $M_{\ge c}(c)$ is generated in degree $0$ its minimal resolution $F$  satisfies $\beta_{k,a}(F)=0$ for $|a| < k$. Thus the first assertion holds for $M_{\ge c}(c)$, because we can obtain $F$ from $F'$ by canceling trivial 
subcomplexes. 

For the second assertion, we note that  in the sheafified complex
$$
0 \leftarrow \cF(c)  \leftarrow \widetilde G_0(c) \leftarrow \widetilde G_1(c) \leftarrow  \ldots 
$$
we have $H^p(\PPn, \widetilde G_k(c))=0 $ for $p\ge 1$. Thus the complex is exact on global sections, 
$$
H^0(\cF(c)) = \coker(H^0(\widetilde G_1(c)) \to H^0(\widetilde G_0(c)))=\coker((G_1)_c \to (G_0)_c) = M_c,
$$
 and the higher cohomology  of $\sF(c)$ vanishes.
A similar argument shows that $R^{p}\pi_{J}\sF(c) = 0$ for $p>0$.
\end{proof}

\begin{proposition}\label{prop:sufficiently positive for pushforward}
If $c$ is sufficiently positive for $M$ then $0$ is sufficiently positive for 
$\Gamma_{\ge 0}({\pi_J}_*( \sF(c))).$
\end{proposition}

\begin{proof}
If $c$ is sufficiently positive for $M$ and $J \subset \{1,\ldots,t\}$ then
$$
\Gamma_{\ge 0}({\pi_J}_*( \sF(c))) = \sum_{a\ge 0, \, a_j=0\, \forall j \notin J} M_{c+a},
$$
where we consider $M_{c+a}$ to be in degree $a$.

  For property 1 of the definition, we note that the subquotient complex of  the free resolution of $M$, consisting of elements that have degrees of the form
$c+a$ with $a\ge 0$ and $a_j=0$ for all $j \notin J$, is a linear free resolution of the module
$\sum_{a\ge 0, \, a_j=0\, \forall j \notin J} M_{c+a}$ over the homogeneous coordinate ring of $\P_{J}$.

For property 2, note that for $a\ge 0$ with $a_j=0\, \forall j \notin J$ we have 
$$
H^0(\PP_J,{\pi_J}_* (\sF(c))(a))=H^0(\PP,\sF(c+a))=M_{c+a}
$$ 
and
$$
H^p(\PP_J,(\pi_J)_* (\sF(c))(a))=H^p(\PP,\sF(c+a))=0,
$$ 
since the direct images $R^p (\pi_J)_*( \sF(c))$  vanish for $p>0$.
Also,
for $I\subset J$ and the further projection $\pi_{I \subset J}: \PP_J \to \PP_{I}$ 
we have 
$$
(\pi_{I \subset J})_*( \pi_J)_*( \sF(c))(a)) = (\pi_I)_*(\sF(c+a))
$$ 
and  
$$
R^p(\pi_{I \subset J})_*(( \pi_J)_* (\sF(c))(a)) = R^p(\pi_I)_*(\sF(c+a))=0
$$ 
for all $a \ge 0$ having $a_j=0$ for $j \notin J$.
\end{proof}

 \begin{notation}[Vectors of Ones] \label{notation:vectors of one}
 Set
\begin{align*}
1^t &:=(1,\ldots,1) \in \ZZ^t \\
0^i1^{t-i}&:=(0 \ldots 0,1,\ldots 1) \in \ZZ^i \times \ZZ^{t-i}.
\end{align*}
We denote the standard basis vectors by $1_i =0^{i-1}10^{t-i}=(0 \ldots,0,1,0\ldots 0) \in \ZZ^t$.
\end{notation}

\begin{notation}[Degree restrictions] Let $M= \sum_{a \in \ZZ^t} M_a$ be a multigraded $S$-module, 
$J \dot\cup I =\{1,\ldots t\}$ complementary subsets of the index set,
 $c=(c_I,c_J)\in \ZZ^{t}$ a multidegree broken up into two parts, and $d\in \ZZ$. We set:
\begin{align*}
&M_{\ge c} = \sum_{a:\,  a \ge c} M_a \, ,\quad M_{c_I,\ge c_J} = \sum_{a=(a_I,a_J) \atop a_I=c_I,\, a_J \ge c_J} M_a \; \hbox{ and }\; M_{>c_I,\ge c_J} = \sum_{ a=(a_I,a_J) \atop a\ge c , \, a_I \not= c_I} M_a\\
&M_{\le c}= \sum_{a:\, a \le c} M_a  \,  \hbox{ and } \, M_{\ge d} =\sum_{a: \, |a| \ge d} M_a. 
\end{align*}

Decomposing $S=S_I\tensor S_J$ accordingly we may regard $M_{c_I,\ge c_J} \cong M_{\ge c}/M_{>c_I,\ge c_J}$ either as an $S$-module or an $S_J$-module. Let $E=E_I \tensor E_J$ be the corresponding decomposition of the exterior algebra and 
$$
\bR_J: \grmod(S_J) \to \lincplx(E_J)
$$ 
the corresponding $\bR$-functor. Note that
$$
\bR(M_{c_I,\ge c_J}) \cong \Hom_K(E_I,K(c_I)) \tensor \bR_J(M_{c_I,\ge c_J}).
$$
is a flat extension of scalars, so 
$\bR(M_{c_I,\ge c_J})$ is exact iff $\bR_J(M_{c_I,\ge c_J})$ is exact.
\end{notation}

\begin{proposition}\label{strands1} Let $\{j\} \cup I = \{1,\ldots,t \} $ 
be a disjoint decomposition with a singleton. 
If $b=(b_j,b_I)$ and $c=(b_j+1,b_I)$ are multidegrees such that $M_{\ge b}(b)$ and $M_{\ge c}(c)$  have  linear resolutions,
then $\bR(M_{\ge b}), \, \bR(M_{\ge c})$ and $\bR(M_{b_j,\ge b_I})$ are acyclic and $M_{b_j,\ge b_I}(b)$ has a linear resolution.
\end{proposition}

\begin{proof} Since $M_{\ge c} \subset M_{\ge b}$ is a submodule with quotient $M_{b_j,\ge b_I}$ we have a short exact sequences of complexes
$$
0 \to   \bR(M_{\ge c}) \to \bR(M_{\ge b}) \to \bR(M_{b_j,\ge b_I}) \to 0
$$
By the Reciprocity Theorem the only non-zero cohomology groups in the long exact sequence are
$$ 
0 \to \HH^{|b|}(\bR(M_{\ge b})) \to \HH^{|b|}(\bR(M_{b_j,\ge b_I}))  \to \HH^{|b|+1}(\bR(M_{\ge c})) \to 0 .
$$
So $M_{b_j,\ge b_I}(b)$ has a linear resolution as well.
\end{proof}

\begin{theorem}\label{cut complex1} Let $c$ be a multidegree  that is sufficiently positive for $M$.
There is a natural \emph{splice map}
$\varphi_{c}: \bR(M_{c}) \to \bR(M_{c+1^{t}})$ that is a surjection to the module
of cycles $P$ of homological degree $|c|+t$ of the complex $\bR(M_{\geq c+1^{t}}).$
Moreover, if 
$c-k\b1$ is sufficiently positive for $M$ for some $k \ge 0$, then $\bR( (M_{\le c})_{\ge |c|-k})$ 
 is the beginning of a projective resolution of $P$.
\end{theorem}

%\david{internally the complex of Theorem~\ref{cut complex1} is called a ``cut complex'' but the term is never defined, so I eliminated the one use of it in the open text. Would it be helpful to introduce the term?}

\begin{proof} We will prove the first statement by induction on the number $t$ of factors in the product $\P$. If $t=1$ then the statement is equivalent to the acyclicity of $\bR(M_{\ge c})$, which follows from the Reciprocity Theorem. 
In this case the splice map $\varphi_{c}$ is simply the differential.

%Now consider the case $t>1$.
%Let $p \in \bR(M_{c+\b1})$ be a cycle in $\bR(M_{\ge c+\b1})$. Since  $\bR(M_{\ge c+\b1})$ is a subcomplex of $\bR(M_{\ge c+0^11^{t-1}})$ it is also a cycle in this complex and by its exactness we find an element $p_1 \in \bR(M_{c+0^11^{t-1}})$ that maps to $p$ under the differential.   In particular $p_1$ is a cycle in
%the complex $\bR(M_{c_1,\ge (c_2+1,\ldots, c_t+1)})$. 
%
%\david{revised version:}
Now consider the case $t>1$. The module $P$  is part of the module
of cycles in the larger complex $\bR(M_{\ge c+0^11^{t-1}})$. By the exactness of this complex, $P$ is
contained in the image of $\bR(M_{c+0^1 1^{t-1}})$; and thus is contained in the image of the
module $P'$ of cycles
of homological degree $|c|+t-1$ in the subquotient complex $\bR(M_{c_1,\ge (c_2+1,\ldots, c_t+1)})$, since
$P$ is congruent to 0 in this subquotient.

By our hypothesis  $M_{c_1,\ge (c_2+1,\ldots, c_t+1)}$ is a truncated section module on a product of $t-1$ factors. By induction, there is a splice map $\varphi'_{c}: \bR(M_{c}) \to \bR(M_{c+0^1 1^{t-1}})$ that is a surjection onto $P'$. Let $\varphi_{c}$ be the composition of this with the component of the differential that maps 
$\bR(M_{c+0^1 1^{t-1}})$ to $\bR(M_{c+1^{t}})$.

By construction, the image of the map $\varphi_{c}$ contains $P$. Thus to prove the first statement
of the Theorem it suffices to show that
the image of $\varphi_{c}$ composes to 0 with the first differential
of $\bR(M_{\geq c+1^{t}}).$
 Because $\varphi_{c}$ is
defined as the composition of 
of $\varphi'_{c}$ with the component of the differential that maps 
$\bR(M_{c+0^1 1^{t-1}})$ to $\bR(M_{c+1^{t}})$, it is at least clear that
the $\varphi_{c}$ composes to 0 with the component of the differential that goes from
$\bR(M_{c+1^{t}})$ to $\bR(M_{c+2^{1}1^{t-1}})$ is 0.

Unraveling the induction used to define $\varphi_{c}$, we see that $\varphi_c$ is the composition
$$
\bR(M_{c}) \to \bR(M_{c+0^{t-1}1^{1}}) \to \ldots \to \bR(M_{c+0^2 1^{t-2}}) \to \bR(M_{c+0^1 1^{t-1}})  \to \bR(M_{c+\b1})
$$
of $t$ maps that are components of  differentials of $\bR(M_{\ge c})$. If we permute the $t$ factors of $\P$, then we would define the map 
$\bR(M_{c}) \to \bR(M_{c+\b1})$ as a different composition. However, since all squares in $\bR(M_{\ge c})$ anti-commute, 
the resulting map would only differ by the signature of the permutation. In particular we see that the composition of 
$\varphi_{c}$ with every component of the first differential of $\bR(M_{\ge c+\b1})$ is zero as required.

\setcounter{MaxMatrixCols}{16}
\begin{figure}
\tiny{
$$
\begin{matrix}
&&&& & & \uparrow & & \uparrow & & \uparrow & & \uparrow & \cr
&&&&&&&&&&&&&\cr
&&&& & & \bR M_{(1,3)}& \rightarrow &  \bR M_{(2,3)}  & \rightarrow &  \bR M_{(3,3)} & \rightarrow & \bR M_{(4,3)} & \rightarrow  \cr
&&&&&&&&&&&&&\cr
&&&& & & \uparrow & & \uparrow & & \uparrow & & \uparrow & \cr
&&&&&&&&&&&&&\cr
&&&& & & \bR M_{(1,2)}& \rightarrow &  \bR M_{(2,2)}  & \rightarrow &  \bR M_{(3,2)} & \rightarrow & \bR M_{(4,2)} & \rightarrow  \cr
&&&&&&&&&&&&&\cr
&&&& & & \uparrow & & \uparrow & & \uparrow & & \uparrow & \cr
&&&&&&&&&&&&&\cr
&&&& & & \bR M_{(1,1)}& \rightarrow &  \bR M_{(2,1)}  & \rightarrow &  \bR M_{(3,1)} & \rightarrow & \bR M_{(4,1)} & \rightarrow  \cr
&&&&&&&&&&&&&\cr
&&&& &\nearrow & &  &  & &  &  &  &  \cr
&&&&&&&&&&&&&\cr
\bR M_{(-2,0)}&\rightarrow&\bR M_{(-1,0)}&\rightarrow& \bR M_{(0,0)}\cr %& & &  &  & &  &  &  &  \cr
&&&&&&&&&&&&&\cr
 &&\uparrow&& \uparrow& & &  &  & &  &  &  &  \cr
 &&&&&&&&&&&&&\cr
 &&\bR M_{(-1,-1)}&\rightarrow& \bR M_{(0,-1)}& & &  &  & &  &  &  &  \cr
 &&&&&&&&&&&&&\cr
 &&&& \uparrow& & &  &  & &  &  &  &  \cr
 &&&&&&&&&&&&&\cr
 &&&& \bR M_{(0,-2)}& & &  &  & &  &  &  &  \cr
\end{matrix}
$$}
\caption{ \small{Nonzero terms of the complex $R((M_{\le 0})_{\ge -2}) \to R(M_{\ge \b1})$ in the case where $t=2$ and where
$M$ is a module for which $b=(-2,-2)$ is sufficiently positive.}}
\end{figure}

We will not give a proof of the second statement of the Theorem here, since the result follows from Theorem~\ref{exactness of T} proven below; in any case, it can be proven by a straightforward diagram chase.
\end{proof}

\begin{corollary}\label{boxed Betti number} Suppose  the multidegree $b$ is sufficiently positive for the multigraded $S$-module $M$ and let $c \ge b+\b1$. If the graded betti number $\beta_{|a|,a}(M)$ is nonzero, then
$0 \le a_i \le n_i$. Moreover,  $\beta_{|n|,n} \not=0$. 
In particular $\depth M_{\ge c} = t$. 
\end{corollary}

\begin{proof}
Let $P$ be the $E$-module of cycles of homological degree $|c|$ in $\bR(M_{\ge c})$.
By the Reciprocity Theorem, the minimal free resolution of $M_{\ge c}(c)$ is the complex $L(P(c)) = L(P)(c)[-c]$, 
so it suffices to examine the Hilbert function of $P$.

The socle of $P$ coincides with the socle of $\bR(M_{ c})$, which is nonzero in degree $c$, while the generators of $P$ have the same degrees as the generators of $\bR(M_{c-\b1})$, that is, $c-n$, and the first statement follows. The last statement follows from the first by the Auslander-Buchsbaum formula.
\end{proof}

Let $M$ be graded $S$-module and let $b$ a sufficiently positive multidegree for $M$. For each $c\ge b+\b1$ we consider 
the module of cycles
$P=\HH^{|c|}(\bR(M_{\ge c})$ as an $E$-module of cohomological degree $|c|$. Let
$$
\Tail_c(M) \to P
$$
be a minimal free resolution of $P$ as an $E$-module. 
By Theorem \ref{cut complex1} the first term of $\Tail_c(M)$ is $\bR(M_{c-\b1})$. By definition, this module has socle in degree $c-1^t$. 

Now suppose that $I$ is the complement of the singleton $\{j\}\subset \{1,\ldots,n \}$ and write
$c+1_{j}$, recalling that $1_j$ was defined in Notation~\ref{notation:vectors of one}. The short exact sequence 
of modules 
$$
0\to M_{\ge c+1_j} \to M_{\ge c} \to M_{c_j,\ge c_I} \to 0
$$ gives rise to a short exact sequence of minimal injective resolutions, from which we deduce
a short exact sequence of modules of cycles:
$$
0 \to \HH^{|c|}(\bR(M_{\ge c})) \to \HH^{|c|}(\bR(M_{c_j,\ge c_I}))  \to \HH^{|c|+1}(\bR(M_{\ge (c+1_j)})) \to 0 .
$$
Corresponding to the left hand map, we get a map $\alpha_{c,j}$ from the minimal free resolution
$\Tail_{c}(M)$ of $\HH^{|c|}(\bR(M_{\ge c}))$ to the minimal free resolution of 
$\HH^{|c|}(\bR(M_{c_j,\ge c_I}))$ whose mapping cone is a free resolution
of $ \HH^{|c|+1}(\bR(M_{\ge (c+1_j)}))$.

\begin{proposition} \label{minimality of mapping cone}
The minimal free resolution
$\Tail_{c+1_{j}}$ of $\HH^{|c|+1}(\bR(M_{\ge (c+1_j)}))$ is isomorphic to the result of
canceling the first map in the mapping cone
of $\alpha_{c,j}$. Thus we may regard
$\Tail_c(M)$ as a quotient complex of 
$\Tail_{c+1_j}(M)$. 
\end{proposition}

\begin{proof} It is immediate from the long exact sequence in homology that 
the homology of the mapping cone is $\HH^{|c|+1}(\bR(M_{\ge (c+1_j)}))$, in homological degree 1.
Since
$\alpha_{c,j}$ is an isomorphism on the terms of homological degree 0, we may cancel these
to obtain a free resolution of $\HH^{|c|+1}(\bR(M_{\ge (c+1_j)}))$. 

It remains to prove minimality. The module  $M_{c_j,\ge c_I}$ is annihilated by $W_{j}$, so the differential of the
complex $\bR(\bR(M_{c_j,\ge c_I})$ does not involve any of the exterior variables in $V_{j}$. Thus
the module $\HH^{|c|}(\bR(M_{c_j,\ge c_I}))$ is free over the tensor factor $E_{j} = \wedge V_{j}$ of $E$,
and its socle has degree $c_{j}$ in the $j$-th component. 

It follows that the free modules in the minimal
free resolution  of $\HH^{|c|}(\bR(M_{c_j,\ge c_I}))$
all have socle with  degree $c_j$ in the $j$-th component. On the other hand it follows from
Theorem~\ref{cut complex1} that
the free modules in the minimal free resolution
 resolution $\Tail_c(M)$ of $\HH^{|c|}(\bR(M_{\ge c}))$
have socle degree in the $j$-th component all $\le c_j-1$. Thus the
mapping cone is minimal as claimed.
\end{proof}

If we choose a sequence of sufficiently positive multi-indices $c^{(i)}$ going to infinity in each component,
then Proposition~\ref{minimality of mapping cone} allows us to give the collection 
$$
\{ \Tail_{c^{(i)}}(M) \mid  c \ge b\}
$$ 
the structure of a directed system $\Tail_{c^{(i+1)}}(M)  \twoheadrightarrow \Tail_{c^{(i)}}(M)$.
We define  the Tate resolution of $\sF= \widetilde M$ as the restricted inverse limit
$$
\bT(\sF) := \left({\lim_{\infty\leftarrow c}}' \Tail_{c}(M)\right)[1-t],
$$
that is, as the complex generated by limit elements that are represented by sequences of elements all
of which are homogeneous of the same degree. The shift $[1-t]$ is necessary to adjust for the shift in the complex
of Theorem~\ref{cut complex1}. 
The complex $\bT(\sF)$ depends, up to isomorphism, only on $\sF= \widetilde M$.

In Corollary~\ref{Tate is Tate} we will show that $\bT(\sF)$ is in fact a Tate resolution in the sense defined in the introduction---that is, all its strands are exact.  For now we prove a weaker property:

\begin{proposition}\label{first properties of Tate} Let $\sF$ be a coherent sheaf on $\P$.
\begin{enumerate}
 \item For each multidegree $a$ the space of homogeneous elements $\bT(\sF)_a$ of multidegree $a$
is finite-dimensional. 
\item For any multi-index $a$ we have 
$$
\bT(\sF(a)) \cong \bT(\sF)(a)[-a].
$$
\item $\bT(\sF)$ exact and uniquely determined by $\sF$ up to isomorphism.
\end{enumerate}
\end{proposition}

\begin{proof}
1. Since each $\Tail_{c}(M)$ is a minimal resolution of a finitely generated module, and the dimension $\dim_K E$ is finite as a $K$-vector space, the space of homogeneous elements $(\Tail_{c }(M))_a$ of degree $a$ is finite dimensional. In the inductive construction of the $\bT(\sF)$ the kernel of the map $\Tail_{c+1_{j}}M \to \Tail_{c} M$ is part of the resolution of a submodule of $\bR(M_{c_{j}, \geq c_{I}})$, so the only multidegrees in which modules in 
the resolution are nonzero have $j$-th coordinate between $c_{j}$ and $c_{j}+n_{j}$. Thus a given degree
can appear in only finitely many kernels, and thus, for every $c$, the dimension $\dim_K (\Tail_{c}(M))_a$ stabilizes as $c \to \infty$. Thus $\bT(\sF)_a$ is a finite dimensional for each multidegree $a$. 

3. Since each $(\Tail_{c}(M))$ is acyclic, this stabilization shows that each $(\Tail_{c}(M))_{a}$ is exact, and it follows that $\bT(\sF)_a$ is exact. 
If the graded $S$-modules $N$ and $M$ represent the same sheaf then $N_{\ge c} \cong M_{\ge c}$ for 
$c \gg 0$. So $\bR(N_{\ge c}) \cong \bR(M_{\ge c})$ and 
$\Tail_c(N) \cong \Tail_c(M)$ as minimal free resolution over a graded ring. Finally $\left({\lim_{\infty\leftarrow c}}' \Tail_{c}(M)\right)[1-t]\cong \left({\lim_{\infty\leftarrow c}}' \Tail_{c}(M)\right)[1-t]$.

2. From Lemma~\ref{TwistShift} it follows that $\bT(\sF(a)) \cong \bT(\sF)(a)[-a]$ for any multi-index $a$. 
\end{proof}

%A more precise statement about the graded pieces is our first main result.
%\david{this is now Corollary~\ref{coho from Tate}. Maybe just refer to it here rather than restating it.}
%
%\begin{theorem}\label{sheaf cohomology} The Tate resolution of a coherent sheaf $\sF$on $\PPn$ has terms
%$$
%\bT(\sF)^d \cong  \bigoplus_{a \in \ZZ^t } \Hom_K(E, H^{d-|a|}(\sF(a)))
%$$
%where we regard $H^{d-|a|}(\sF(a))$ as a vector space in degree $a$
%\end{theorem}
%
%The proof of this result relies on the Beilinson monad, and will be given in the next section.

% \section{Beilinson Monads}
\section{Beilinson Monads}\label{monads}

To prove the more precise statement about the graded pieces of the Tate resolution
given in Theorem~\ref{Tate existence}
we will use  Beilinson monads.
Recall that a \emph{Beilinson monad} for a coherent sheaf $\sF$ on $\P$ is a finite complex $B$ whose only homology is $\sF$ in homological degree 0, and whose
terms that are direct sums of the sheaves $U^{a}$. It is \emph{minimal} if, under the identifications $\Hom_{\P}(U^{a}, U^{b}) = E_{b-a}$, all the nonzero components of the 
morphisms of $B$ have nonzero degree. The following result is well-known:

\begin{theorem}\label{Terms of the Beilinson monad}  If $B$ is a minimal Beilinson monad for a coherent sheaf $\sF$ on $\P$, then the $d$-th term is
$$
B^d =  \bigoplus_{a \in \ZZ^t} H^{d-|a|}(\sF(a)) \tensor U^{-a} 
$$
where (in this formula) we regard $H^{d-|a|}(\sF(a))$ as a vector space in degree $0$.
In particular, $B^{d}$ can be nonzero only in the range $-|n|\leq d \leq |n|$.
\end{theorem}

\begin{proof} Write $B^d = \oplus_{0\le a \le n}  B_{a}^d \tensor U^{-a}$ for vector spaces $B_a^d$ in degree zero. The collection $\{U^{-a} | -n \le a \le 0 \}$ forms a strong full exceptional collection
for the derived category $D^b(\PP)$ of coherent sheaves on $\PP=\PPn$ (see \cite[Def. 8.31]{huybrechts}) which is right orthogonal to the strong full exceptional collection
$\{\sO(c) | 0 \le c \le n \}$ in the sense that, for $c$ in this range,
$$
H^p {\RHom(\sO(c),U^{-a})} = H^{p}( U^{-a}(-c)) =
\begin{cases}
 K & \hbox{ if $c=-a$ and $p=-|a|$, }  \cr
0 & \hbox{ otherwise.}
\end{cases}
$$
Hence the complex $\RHom(\sO(a),\bU(\sF))$ is the sequence 
$$
 \ldots \to B_a^d \to B_a^{d+1} \to \ldots 
$$
of vector spaces, where $B_a^d$ sits in cohomological position $d-|a|$, and all maps are zero since $B$ is minimal. 
Since $B$ is quasi isomorphic to $\sF$ this complex coincides with
$\RHom(\sO(a),\sF)$ which is a sequence of vector spaces with terms
$$
 0 \to H^0(\sF(a)) \to H^1(\sF(a)) \to \ldots \to H^{|n|}(\sF(a))  \to 0
$$ 
in cohomological position $0, \ldots, |n|$ and zero differential. Thus
$$
B_a^d=H^{d-|a|} \sF(b)
$$
as desired. Note that this group is possibly nonzero only if $-|n| \le d \le |n|$, because $0 \le d-|a|  \le |n|$ and $0 \le a \le n$.
\end{proof}

% Recall that $\omega_E=\Lambda W= \Hom_K(E,K) $ is the injective hull of $K = E/ \gm_E$. It plays the role of a dualizing module.
% Each term in the Tate resolution can be written a a direct sum of $\omega_E(a)$ for various multidegrees $a$.

%Let $U_k = \ker(H^0(\PP^{n_k},\sO(1))\tensor \sO \to \sO(1))$ denote  the rank $n_k$ tautological subbundle on  $\PP^{n_k}$. As in the introduction, for $a\in \ZZ^{t}$, we denote by
%$$
%U^a = \boxtimes_{k=1}^t \Lambda^{a_k} U_k=(\pi_1)^*\Lambda^{a_1}U_1 \otimes \cdots \otimes (\pi_t)^*\Lambda^{a_t} U_t
%$$
%the tensor product of the pullbacks to $\PP$ of exterior powers of the $U_k$. 
%Of course $U^a$ is nonzero if only if $0 \le a \le n$. 
%
%We now define an additive functor from the category of finitely generated free graded right $E$-modules 
%to the category of finite direct sums of the sheaves $U^{a}$.  The functor is defined on objects by:
%$$
%\bU\colon \omega_E(a) \mapsto U^a.
%$$
%A morphism in the category of free right $E$-modules
%may be given by a matrix, acting on the left, so it suffices to define the functor
%on morphisms of the form $e\wedge -: E(a)\to E(b)$, where  $e \in E_{b-a} = \otimes_{k=1}^t \Lambda^{a_k-b_k}V_k$. We define
%$U(e\wedge -) :=  \neg e$, where $\neg e$ is contraction with $e$, as in the following diagram:
%$$
%\xymatrix{
%\omega_E(a)  \ar[d]_{e\wedge-}& \mapsto & U^a \ar[d]^{\neg e} & \hookrightarrow & \otimes_{k=1}^t \Lambda^{a_k}W_k\otimes_K \sO \ar[d] ^{\neg e}\cr
%\omega_E(b) &\mapsto & U^b & \hookrightarrow & \otimes_{k=1}^t \Lambda^{b_k}W_k \otimes_K \sO \cr },
%$$

We can now show that $\bU(\sF):= \bU(\bT(\sF))$ is a Beilinson monad for $\sF$.
Because each $\bT(\sF)_{a}$ is finite dimensional, the result of applying $\bU$ to the Tate resolution of $\sF$ is at least a bounded complex.

\begin{theorem}\label{Beilinson monad from Tate}
$\bU(\sF)$ is a minimal Beilinson monad for the sheaf $\sF$.
\end{theorem}

\begin{proof} As in the case of a single projective space, the result follows by applying the functor $\bL$ to a bounded part of the Tate resolution.

We first consider the complexes $\bL(\omega_E(a))$. For an $E$-module $Q=\sum_a Q_a$ we set $Q_{\ge 1^t}:=\sum_{a \ge 1^t} Q_a$. If we write $Q_a= B_a(-a)$ with $B_a$ 
a vector space in degree $0$, then the sheafification $\widetilde \bL(Q_{\ge 1^t})$ of the complex $\bL(Q_{\ge 1^t})$  has the form
$$
\cdots \to \oplus_{j} B_{1^{t}+1_{j}}\otimes \sO((-1)^{t}-1_{j})
\to B_{1^t}\otimes \sO((-1)^{t})
\to 0
$$
It is more convenient to index the complex cohomologically, so that this term becomes $\widetilde \bL^{-t}(Q_{\ge 1^t})$. 
For example 
$$
\widetilde \bL((\omega_E)_{\ge \b1}):  \ldots \to (\otimes_{k=1}^t W_k) \tensor \sO(-1,\ldots,-1) \to 0
$$
is the tensor product of $t$ truncated Koszul complexes
$$
(\pi_k)^*(0 \to \Lambda^{n_k+1} W_k \tensor \sO(-n_k-1) \to \ldots \to W_k \tensor \sO(-1) \to 0).
$$
Since each of these complexes is a resolution of $\sO$ we see that
$$
\HH^{-t}(\widetilde \bL((\omega_E)_{\ge \b1})) = \sO
$$
and all other cohomology groups are zero. Similarly we find that the rightmost term of
$
\widetilde \bL((\omega_E(a))_{\ge \b1})
$
is 
$(\otimes_{k=1}^t \Lambda^{a_k+1} W_k) \tensor \sO(-1,\ldots,-1)$
and the whole complex is a tensor product of possibly truncated Koszul complexes. Thus the complexes 
$
\widetilde \bL((\omega_E(a))_{\ge \b1})
$
are acyclic with only cohomology
$\HH^{-t}(\widetilde \bL((\omega_E(a))_{\ge \b1})= U^a$, which is
nonzero if and only if $0 \le a \le n$. 

For any sheaf $\sF$ we can form the double complex
$$
\widetilde \bL((\bT(\sF))_{\ge 1^t})[t].
$$
What we have proven shows that the  homology of this double complex with respect to the differential coming from the functor $\bL$, which we think of as the horizontal homology, 
is the complex $\bU(\sF)$. 

Let $M$ be a graded $S$-module whose sheafification is $\sF$.  The complex $\bU(\sF)$ 
depends only on the summands $\omega_{E}(a)$ of $\bT(\sF)$ 
with $0\leq a\leq n$. 
Choose a multidegree $c$ that is sufficiently positive for $M$ and large enough such that every such summand 
is contained in $\Tail_{c}(M)$.  
 With this choice, $\widetilde \bL((\Tail_{c}(M))_{\le \b1})$ is a bounded double complex whose horizontal cohomology is
$\bU(\sF)[-t]$. 

Recall that $\Tail_{c}(M)$ is, by definition,
 a resolution of the module of cycles $P=\HH^{|c|}(\bR(M_{\ge c}))$, which has cohomological degree $|c|-t$.  Thus the  vertical homology
 of $\widetilde \bL((\Tail_{ c }(M))_{\le \b1})$ is the complex 
 $$
\widetilde \bL(P)[-t+|c|]
$$ 
whose only homology is $\HH^{-t}(\widetilde \bL(P)[-t+|c|]) \cong \widetilde M_{\ge c} \cong \sF$.

It follows that the only homology of $\bU(\sF)$ is $\HH^0(\bU(\sF))\cong \sF$. More precisely, the total complex of $\widetilde \bL((\Tail_{c}(M))_{\ge 1^t})[t]$ is quasi isomorphic
to both $\bU(\sF)$ and $\widetilde \bL(P )[-|c|]$ via the natural maps
$$
\begin{diagram}
&&\widetilde \bL((\Tail_{c}(M))_{\ge 1^t})[t] \\
&\ldTo&&\rdTo\\
\widetilde \bL(P )[-|c|]&& &&\bU(\sF).
\end{diagram}
$$
\end{proof}

\begin{proof}[Proof of Theorem~\ref{Tate existence}]
The complex $\bT(\sF)$ is minimal by construction, and we will see that it is a Tate resolution in Corollary~\ref{Tate is Tate} below.
Applying Theorem~\ref{Beilinson monad from Tate} and Theorem~\ref{Terms of the Beilinson monad} to the sheaf $\sF(b)$ we see that, if we write
$\bT(\sF(b))^{d} = \oplus_{a}B^{d}_{a}\otimes \omega(-a)$
then for $-n\leq a \leq 0$ we must have $B^{d}_{a} = H^{d-|a|}\sF(a+b)$.
Moreover, by Proposition~\ref{first properties of Tate}, 
$$
\bT(\sF) = \bT(\sF(b))(-b)[-b]
$$
for every multi-index $b$.  It follows that $\bU(\bT(\sF)(b)[b])$ is a Beilinson monad for $\sF(b)$. $\bT(\sF)$ is uniquely determined by $\sF$ up to isomorphism by Proposition \ref{first properties of Tate}. It is also uniquely determined up to isomorphism by any of its Beilinson monads $\bU(\bT(\sF)(b)[b])$ by Corollary \ref{reconstruction} below.
\end{proof}
%\section{Exactness properties of the Tate resolution}
\section{Exactness properties of the Tate resolution}\label{exactness properties}

We will next establish the exactness of the strands and other subquotient complexes of $\bT(\sF)$, showing that it
is indeed a Tate resolution as defined in the Introduction. We begin by establishing
notation. First, we can restrict the notions already defined to any nonempty subset
$J$ of the 
indices  $\{1,\dots,t\}$:
We set 
 $
 W_J := \oplus_{j \in J} W_j,
 $
and $S_{J} := \Sym W_{j}$ We use similar notation for $V$ and $E$.
As in the Introduction we write
$
\pi_J: \PP \to \PP_J= \prod_{j \in J} \PP^{n_j}.$

We denote by $\omega_{J} = E_{J}(-\dim V_{J})$ the $E_{J}$-injective hull of
$K$, and we write $\bU_{J}$ for the functor
 whose value on $\omega_{J}(a)$ is $\otimes_{j\in J} \wedge^{a_{j}}U_{j}$,
 analogous to $\bU$.

%We denote by 
%$$
%e_J \in \Lambda^{\dim V_J} V_J \tensor K  \subset E_J \tensor K \subset E
%$$
%the generator of the socle of $E_J$,  regarded as an element of $E$, and by
%$\gm_J =\Ann_E e_J$ its annihilator.
%Thus $\gm_J$ is the ideal generated by $V_J$. 

%\begin{definition} Let $T$ be an complex of graded free $E$-module,with terms
%$$
%T^d = \sum_{a \in \ZZ^t} B^d_a \tensor \omega_E(-a)
%$$
%with vector spaces $B_a^d$ in degree $0$, is \textbf{ locally finite}, if
%for each $a$ the vector space 
%$
%\sum_{d\in \ZZ} B^d_a
%$
%is finite dimensional.
%\end{definition}

%Let $J \cup I =  \{1,\dots,t\}$ be a decomposition into disjoint subsets,
%and let $T$ be a complex of graded free $E$-module having terms
%$$
%T^d = \sum_{a \in \ZZ^t} B^d_a \tensor \omega_E(-a)
%$$
%with vector spaces $B_a^d$ in degree $0$. We write
%$T_I$ for the subquotient complex with terms
%$$
%T_I^d = \sum_{a=(a_I,a_J)\in \ZZ^t \atop a_I=0} B^d_{a} \tensor \omega_{E_J}(-a_J).
%$$
%Thus the differential of $T_{I}$ is given by matrices in with entries in $E_J$, and 
%$$
%T_{I}\otimes \omega_{E_{I}}
%$$
%is the $I$-th strand of $T$. (In the notation of
% \ref{strands etc}, below, we have $T_{I}\otimes \omega_{E_{I}} = \Tate 0 \o {I} \o$.)
%\end{notation}

%\begin{remark}
%A candidate for $R {\pi_J}_*(\bU(T)(c))$ is $\bU_J((T(c)[c])_{I})$, but this holds only under further hypotheses on $T$.
%\end{remark}

\begin{notation} [Strands, quadrant complexes, and region complexes]\label{strands etc}
Let $T$ be a locally finite complex of graded free $E$-modules with terms 
$
T^d = \sum_{a \in \ZZ^t} B^d_a \tensor \omega_E(-a).
$
For $c \in \ZZ^t$ and three disjoint subsets $I,J,K \subset \{1,\ldots,t\}$ we denote by $\Tate c I J K$ the subquotient complex
with terms
$$
\Tate c I J K ^d= \sum_{a \in \ZZ 
          \atop {a_i < c_i \hbox{ \scriptsize for } i \in I
         \atop {a_i  = c_i \hbox{ \scriptsize for } i \in J
         \atop a_i  \ge c_i \hbox{ \scriptsize for } i \in K}}} 
          B^d_a \tensor \omega_E(-a),
$$
and we call this a \emph{region complex} of $T$.  A \emph{strand} of $T$, which was
defined in the Introduction, may be viewed as a region complex of the special form
$\Tate c \o J \o$ where $J \subsetneq \{1,\ldots,t\}$ is a proper subset.  Note that 
$T$ itself is the strand corresponding to $J=\emptyset$.
% We write
%$T_I$ for the subquotient complex with terms
%$$
%T_I^d = \sum_{a=(a_I,a_J)\in \ZZ^t \atop a_I=0} B^d_{a} \tensor \omega_{E_J}(-a_J)
%$$
%Thus the differential of $T_{I}$ is given by matrices in with entries in $E_J$, and 
%$$
%T_{I}\otimes \omega_{E_{I}}
%$$
%is the $I$-th strand of $T$. (In the notation of
% \ref{strands etc}, below, we have $T_{I}\otimes \omega_{E_{I}} = \Tate 0 \o {I} \o$.)

If $I\cup J \cup K = \{1,\ldots,t\}$ we call $\Tate c I J K $ a \emph{quadrant complex}.
A region complex which is not a quadrant complex is called a \emph{proper region complex}.
If $T=\bT(\sF)$ is the Tate resolution of a sheaf we will see, that
any proper region complex $\Tate c I J K$ is exact.

To simplify the notation we sometimes write the quadrant complexes as
$$
qT_{c,I}:= \Tate c I \o J
$$
where $J$ is the complement of $I$, and 
 $T_{\ge c} = \Tate c \o \o {\{1,\ldots,t\}}$ and $T_{<c} = \Tate c {\{1,\ldots,t\}} \o \o $ for the first and last quadrant complex.
 \end{notation}
 
Inclusions of regions give various short exact sequences of complexes. For $i \notin I \cup J \cup K$ we have exact sequences
\begin{align*}
0 &\to \Tate c I J {K\cup\{i\}} \to \Tate c I J K \to \Tate c {I\cup\{i\}} J K \to 0,\\
0 &\to \Tate {{c+1_i}} I J {K\cup\{i\}} \to \Tate c I J {K\cup\{i\}} \to \Tate c I {J\cup\{i\}} K \to 0,\\
0 &\to \Tate c I {J\cup\{i\}} {K} \to \Tate {{c+1_i}} {I\cup\{i\}} J K \to \Tate c {I\cup\{i\}} J K \to 0.
\end{align*}
%$$
%0 \to \Tate {{c-1_{i}}} {I\setminus \{i\}} {\{i\}} \o 
%\to \Tate c {I} \o \o
% \to \Tate {{c-1_{i}}} {I} \o \o \to 0
%$$

\begin{notation}[Corner complexes]\label{corner}
We define the map 
$$
\varphi_{c}: T_{<c}[-t] \to T_{\ge c}
$$
to be the composition
$$
 T_{<c}[-t]=qT_{c,\{1,\dots,t\}}[-t]  \to \ldots \to qT_{c,\{1,\dots k\}}[-k] \to \ldots \to qT_{c,\o}= aT_{\ge c}
 $$
of $t$  morphisms deduced from the exact sequences
$$
0 \to qT_{c,\{1,\dots k-1\}} \to \Tate c {\{1,\ldots,k-1\}} \o {\{k+1,\ldots, t\}} \to qT_{c,\{1,\dots k\}} \to 0.
$$ 
This is a morphism of complexes since all  maps $T_{c,k}[-k] \to  T_{c,k-1}[-k+1]$ are morphisms of complexes.
We define the \emph{corner complex} $T_{\Rsh \kern -1pt c}$ as the mapping cone of $\varphi_{c}$. 
\end{notation}

\begin{theorem}\label{exactness of T} Let $T$ be a locally finite complex of free $E$-modules. The following statements are equivalent 
\begin{enumerate}
\item Every strand of $T$ is exact.
\item Every proper region complex of $T$ is exact.
\item Every corner complex $T_{\Rsh \kern -1pt c}$ is exact.
\item The corner complexes $T_{\Rsh \kern -1pt c}$ are exact for every sufficiently large $c$.
\item The proper region complexes $\Tate c I \o \o$ are exact for every sufficiently large $c$.

\end{enumerate}

 \end{theorem}

\begin{proof} By the definition of local finiteness, the complex of vectorspaces $T_{b}$ consisting of homogeneous elements of  degree $b$  is a finite complex of finite dimensional vector spaces. Note that $T$ is exact if and only if $T_{b}$ is exact 
for every internal degree $b$.

$1. \Rightarrow 2.$ For a fixed $c$, we must prove the exactness of $\Tate {{c}} I J K _{b}$ for all $b$.  We note that, by hypothesis, $\Tate {{c}} \o J \o $ is exact for every proper subset $J$.  Let us consider the case $I=\{k\}$ and $K=\o$. If $b_k\ll c_k$ sufficiently large, then we have $\Tate {{c}}  {\{k\}} J \o _{b} =\Tate {{c}} \o J \o _b$, so the 
exactness of $\Tate {{c}}  {\{k\}} J\o_{b}$ follows from the exactness of $\Tate {{c}} \o J \o $. For smaller $c_k$ we use
descending induction on $c_k$ together with the exactness of the strands $\Tate c \o {\{k\}} \o $ and exact sequences
$$ 
0 \to \Tate c \o {J\cup{\{k\}}} \o _{b} \to \Tate {{c+1_k}}  {\{k\}} J \o _{b}  \to \Tate {{c}}  {\{k\}} J \o  _{b} \to 0.
$$
Note that the term on the left is also exact by hypothesis so long as $J\cup \{k\}$ is still a proper subset, yielding the exactness of the term on the right.

The exactness of complexes $\Tate {{c}}  \o J {\{k\}}$ can be argued similarly, using the exact sequence
$$ 
0 \to \Tate c \o {J\cup\{k\}} \o \to  \Tate {{c}} \o  J {\{k\}}    \to  \Tate {{c+1_k}}   \o J {\{k\}}  \to 0
$$
and an ascending induction on $c_k$. 

The general case follows by an induction on the size of $I \cup J \cup K$
using the exact sequences above Notation~\ref{corner}. Note that this induction stops with $\#(I\cup J \cup K)=t-1$ 
because the complex $\Tate c \o {\{1,\ldots,n\}} \o$ is a bounded complex of $E$-modules, and thus is never exact unless $T$ is a split exact complex.

%$1. \Rightarrow 2.$ For a fixed $c$, we must prove the exactness of $\Tate {{c}} I J K _{b}$ for all $b$.  We first consider the case $I=\{k\}$. If $b_k\ll c_k$ sufficiently large, then we have $\Tate {{c}}  {\{k\}} \o \o _{b} =T_{b}$, so the 
%exactness of $\Tate {{c}}  {\{k\}} \o \o_{b}$ follows from the exactness of $T_b$. For smaller $c_k$ we use
%descending  induction on $c_k$ together with the exactness of the strands $\Tate c \o {\{k\}} \o $ and exact sequences
%$$ 
%0 \to \Tate c \o {\{k\}} \o _{b} \to \Tate {{c+1_k}}  {\{k\}} \o \o _{b}  \to \Tate {{c}}  {\{k\}} \o \o  _{b} \to 0.
%$$
%
%The exactness of complexes $\Tate {{c}}  \o \o {\{k\}}$ can be argued similarly, using the exact sequence
%$$ 
%0 \to \Tate c \o {\{k\}} \o \to  \Tate {{c}} \o  \o {\{k\}}    \to  \Tate {{c+1_k}}   \o \o {\{k\}}  \to 0
%$$
%and an ascending induction on $c_k$. 
%
%The general case follows by an induction on the size of $I \cup J \cup K$
%using all  exact sequences above. Note that this induction stops with $\#(I\cup J \cup K)=t-1$ 
%because the complex $\Tate c \o {\{1,\ldots,n\}} \o$ is a bounded complex of $E$ modules, and thus is never exact unless $T$ is a split exact complex. 

$2.  \Rightarrow 3.$ 
By Notation~\ref{corner}, we have that $T_{\Rsh \kern -1pt c}$ is the mapping cone of the map $\phi_c\colon T_{<c}[-t]\to T_{\ge c}$.  Note that, by definition $\phi_c$ is the composition of morphisms $\phi_{c,k}$ arising via the exact sequences
\[
0 \to qT_{c,\{1,\dots k-1\}} \to \Tate c {\{1,\ldots,k-1\}} \o {\{k+1,\ldots, t\}} \to qT_{c,\{1,\dots k\}} \to 0.
\]
In fact, the middle term is the mapping cone of $\phi_{c,k}\colon qT_{c,\{1,\dots k\}}[-1]\to qT_{c,\{1,\dots k-1\}}$.  By hypothesis, the middle term is exact and hence, by the induced long exact sequence, we see that $\phi_{c,k}$ is quasi-isomorphism.  Since a composition of quasi-isomorphisms is a quasi-isomorphism, it follows that $\phi_c$ is a quasi-isomorphism.  By the induced long exact sequence for the mapping cone of $\phi_c$, we then conclude that the mapping cone of $\phi_c$, which equals $T_{\Rsh \kern -1pt c}$, is exact.

$3.  \Rightarrow 4.$ is trivial.

$4.  \Rightarrow 5.$ The complex of vector spaces $(\Tate c I \o \o)_{a}$ obtained by fixing an internal degree $a$
in the region complex $\Tate c I \o \o$  is the same as the degree $a$ part
$(T_{\Rsh \kern -1pt b})_{a}$ in the corner complex $T_{\Rsh \kern -1pt b}$
for $b_{I} = c_{I}$ and $b_{j}>-a_{j}+n$ for some $j\in I'$. 
This is because, for such $b$, the free modules in the the complex
$T_{\geq b}$ do not contain elements of degree $a$.

$5.  \Rightarrow 1.$: 
For any $I\subsetneq \{1,\dots, t\}$ and $c$, the complex $\Tate c {I} \o \o$ is exact for sufficiently large $c$ by hypothesis. By Lemma~\ref{move to front} below, the subcomplex $\Tate {{c-1_{i}}} {I\setminus \{i\}} {\{i\}} \o$ is  exact
for each $i\in I$. From the exact 
sequence
$$
0 \to \Tate {{c-1_{i}}} {I\setminus \{i\}} {\{i\}} \o 
\to \Tate c {I} \o \o
 \to \Tate {{c-1_{i}}} {I} \o \o \to 0
$$
 we see that $\Tate {{c-1_{i}}} {I\cup\{i\}} \o \o$ is exact too. Descending induction now shows that
 $\Tate c {I} \o \o$ is exact for all $c$. 
 
 By Lemma~\ref{move to front} , 
 $\Tate c {I} J \o$ is exact for all $c,I,J$, so long as $I\cup J\subsetneq \{1,\dots, t\}$.  In particular the strands of $T$ are all exact.
\end{proof}

\begin{lemma}\label{move to front}
If $\Tate c  {I\cup J} K \o$ is exact, then so is its subcomplex
$\Tate {{c-1_{J}}} I {J\cup K} \o$. 
\end{lemma}

\begin{proof}
To prove the exactness at $\Tate {{c-1_{J}}} I {J\cup K} \o^{d}$
we decompose $\Tate c {I\cup J} K \o^{d-1}$
as a direct sum of graded free $E$-modules
${\Tate c {I\cup J} K \o'}^{d-1}\oplus \Tate {{c-1_{J}}} I {J\cup K} \o^{d-1}$. Since $\Tate c {I\cup J} K \o$ is 
assumed exact,
the module of cycles $Q\subset \Tate {{c-1_{J}}} I {J\cup K} \o ^{d}$ is equal to
 the sum of the boundaries in $\Tate {{c-1_{J}}} I {J\cup K} \o ^{d}$
and the intersection $P$ of the image of ${\Tate {{c-1_{J}}} I {J\cup K} \o'}^{d-1}$ with $\Tate {{c-1_{J}}} I {J\cup K} \o ^{d}$.
It thus suffices to show that $P$ is contained in the maximal ideal of $E$ times $Q$.

The differential in $\Tate {{c-1_{J}}} I {J\cup K} \o$ involves only the variables of $E_{(J\cup K)'}$, and thus the module of cycles of $\Tate {{c-1_{J}}} I {J\cup K} \o$ in $\Tate {{c-1_{J}}} I {J\cup K} \o ^{d}$ is generated by linear combinations
of the free generators of $\Tate {{c-1_{J}}} I {J\cup K} \o ^{d}$ with coefficients in $E_{L}$, where $L$ is the complement
of $J\cup K$.

On the other hand, it is clear from the form of the differential of $\Tate c {I\cup J} K \o$ that
$P\subset \gm_{E_{J}} \Tate {{c-1_{J}}} I {J\cup K} \o ^{d}$, where
$\gm_{E_{J}}$ denotes the maximal ideal of $E_{J}$. Thus no element of $P$ can be a minimal generator of $Q$.
\end{proof}

\begin{corollary}\label{Tate is Tate} For any coherent sheaf $\sF$ on $\PP$, the complex $\bT(\sF)$ satisfies
 the equivalent conditions of Theorem \ref{exactness of T}. In particular all strands $\bT(\sF)$ are exact,
 and hence $\bT(\sF)$ is a Tate resolution.
\end{corollary}

\begin{proof}
By construction $T$ is a  subcomplex of an inverse limit of the
acyclic complexes $\Tail_{c}(M)$, which are defined for sufficiently positive $c$.
Furthermore $\cornerT c $ coincides with the exact complex obtained as the mapping cone of
$$
\Tail_{c}(M)[-t] \to \bR( \Gamma_{\ge c}(\sF)),
$$ 
so $T$ satisfies condition $4.$ of Theorem~\ref{exactness of T}. 
\end{proof}
The following proposition implies Corollary~\ref{strands}.

%\begin{corollary}
% If $\sF$ is a coherent sheaf on $\PP$, and $T = \bT(\sF)$ is it's Tate resolution, then for each proper subset $J\subset \{1,\dots,n\}$ with complement $I$ and for all degrees $c$ the strand
% $\Tate c \o I \o$ is exact and the complex $\bU_{J}((\Tate c \o I \o(c)[c])_{J})$ \david{has this notation ever
% been defined? The usage is different in Prop~\ref{direct image}} is a Beilinson representative of 
% $R{\pi_J}_{*}(\sF(c))$.
%\end{corollary}

\begin{proposition}\label{direct image} Let $T$ be a locally finite minimal complex of graded free $E$-modules.
%then $\bU(T)$ is the Beilinson representative of a bounded complex of coherent sheaves. 
Let $I\cup J = \{1,\dots, t\}$ be a decomposition into disjoint sets, and let
$T_I$ denote the  complex of $E_{J}$-modules such that
$$
T_{I}\otimes \omega_{E_{I}}
$$
is the $I$-th strand $\Tate 0 \o {I} \o$ of $T$.  The complex $\bU_J(T_I)$
is a Beilinson monad for $R {\pi_J}_*(\bU(T))$ in $D^b(\PP_J)$.
%  has Beilinson monad .
\end{proposition}

\begin{proof}
%Recall that $\bU(\omega_E(a)) = U^a$ and hence $\bU_J(\omega_{E_J}(a_J)) = U_J^{a_J}$.
%Now,
We have
$$
R^p{\pi_J}_* U^a = \begin{cases}
U_J^{a_J} & \hbox{ if $a_I =0$  and } p=0, \cr
0 &\hbox{ otherwise. }\cr
\end{cases}
$$
Since, in particular, $\bU(T)$ is ${\pi_J}_*$-acyclic,  $R{\pi_J}_*(\bU(T))$ is represented by the complex 
$\pi_{J*}\bU(T) = \bU_J(T_I)$. 
%\frank{ fix notation?}
\end{proof}

Proposition~\ref{direct image} gives an interpretation of the 0-strand of a Tate resolution. Theorem~\ref{twist and shift}, which we now prove, provides a similar interpretation for every strand.
%
%\begin{theorem}[Twist and shift formula]\label{twist and shift}
%Let $T$ be a Tate resolution. Then
%$
%\bU(T(c))[c]
%$
%and
%$\bU(T)(c)$
%are quasi-isomorphic. In particular, for any non-empty subset $J \subset \{1,\ldots t \}$ with complement $I$,
% the $I$-th strand of $T$ through $c$ computes the direct image
%$R{\pi_J}_*\bigr( \bU(T)(c)\bigl)$ along $\pi_J: \PP \to \PP_J$ by the formula of Proposition~\ref{direct image}. 
%\end{theorem}
%
\begin{proof}[Proof of Theorem~\ref{twist and shift}]
We may assume $c \ge 0$. Let $b = c+1^t$ and let $P^\bullet$ the image of the lower quadrant complex $T_{\le  c} $ in the upper quadrant complex $T_{\ge b}$ along the corner map. 
Then by the argument
of Theorem \ref{Beilinson monad from Tate}, $\bU(T)$ and $\widetilde{\bL}(P^\bullet)$ are quasi-isomorphic, and similarly $\bU\big(T[c](c)\big)$ and $\widetilde{\bL}(P^\bullet[c](c))$
are quasi-isomorphic. By Lemma~\ref{TwistShift} we have $\bL(P^\bullet[c](c))=\bL(P^\bullet)(c)$, and the first result follows.

The second statement follows  from Proposition \ref{direct image}.
\end{proof}

\begin{example}\label{special corner complexes}

For $\sF=U^a$ and $T=\bT(U^a)$ the corner complex $\cornerT {-a}$ of $T$ at $-a$ has terms
$$
 \ldots \to  \tensor_{j=1}^t \Lambda^{n_j-a_j+1}  V_j \tensor \omega_E(a) \to \omega_E(a) 
\to \oplus_{j=1}^t \Lambda^{a_j+1} W_j \tensor \omega_E(a) \to \ldots
$$
with cohomological indexing such that $(\cornerT {-a} )^0 = \omega_E(a)$. 
Notice that the socle of $\tensor_{j=1}^t \Lambda^{n_j-a_j+1}  W_j \tensor \omega_E(a)$ sits in a single degree $n+\b1$, while the socle of the right hand side sits in several degrees
$a-a_j-1_j$. 
The functor $\bU$ takes all but the middle term to zero. Thus
$\bU(\cornerT {-a}) = U^a$.

The corner complex $\cornerT {1^t}$  has terms
$$
\cdots \to 
\oplus_{j=1}^t \Lambda^{n_{j}-a_j+1} V_j \tensor \omega_E(a) 
\to 
\omega_E(a) 
\to 
\tensor_{j=1}^t \Lambda^{a_j+1}  W_j \tensor \omega_E(a)
\to \cdots
$$
and  satisfies $\bU(\cornerT {1^t}) = U^a[1]$  as well.

\end{example}

\begin{definition} Let $\ell_1 \le \ell_2$ be integers.  A locally finite free complex $T$ of free $E$-modules has \emph{finite amplitude} $[\ell_1,\ell_2]$, if  for all $a \in \ZZ^t$ and $d \in \ZZ$, when we write 
$$
T^d = \sum_a B_{a}^d \tensor \omega_E(-a)
$$
we have that $B^d_a \not=0$
only if $\ell_1 \le d-|a| \le \ell_2$.
\end{definition}

Note that if $T$ has finite amplitude $[\ell_1,\ell_2]$, then $T(c)[c]$ has finite amplitude $[\ell_1,\ell_2]$ as well, while $T[c]$ or $T(c)$ have a different amplitude. The Tate resolution of a coherent sheaf has amplitude $[0, |n| ]$ by Theorem~\ref{Tate existence}.

\begin{corollary}
 Every Tate resolution has a finite amplitude.
\end{corollary}

\begin{proof}  
$\bU(T(c)[c])$ and $\bU(T)(c)$
represent the same object in $D^b(\PP)$ for every $c$. In particular the homology sheaves $H^k(\bU(T(c)[c]))$ occur only in finitely many places, say from $k_1$ to $k_2$ independent of $c$. Thus the hypercohomology $\HHH^{i}(\bU(T(c)[c])$ is nonzero only if $i\in [\ell_1,\ell_2]=[k_1,k_2+|n|]$, and the terms of $T$ have the form
$T^{d} = \oplus_{a}\HHH^{i-|a|}(\bU(T)(c))\otimes \omega_{E}(-a)$ by Theorem~\ref{Terms of the Beilinson monad}.\end{proof}

\begin{remark} If $T$ is a Tate resolution with amplitude $[\ell_1,\ell_2]$, then its Beilinson representative $\bU(T)$ can only have non-zero terms in cohomological degrees
$[\ell_1-|n|,\ell_2]$
\end{remark}

Maclagan and Smith have defined a notion of regularity for a multigraded module over the Cox ring of a smooth toric variety~\cite[Definition~1.1]{maclagan-smith}.  Using the corner complex and reciprocity, we  compare this to the notion of sufficiently positive from Defintion~\ref{defn suff positive}.  Although these notions can diverge (see Remark~\ref{rmk:divergent examples}), they are closely related.  We make use of Notation~\ref{notation:vectors of one} in the following.
\begin{proposition}\label{prop: comparing multigraded regularity and suff positive}
Let $M$ be a finitely generated, multigraded $S$-module such that no associated prime of $M$ contains the irrelevant ideal of $S$.  If $b$ is sufficiently positive for $M$, then $b$ lies in the multigraded regularity of $M$.  Conversely, if $b$ lies in the multigraded regularity of $M$, then for any $i$, $b+1^t-1_i$ is sufficiently positive for $M$.
\end{proposition}
\begin{proof}
Let $\mathfrak b$ be the irrelevant ideal of $S$.  Then $H^0_{\mathfrak b}M = 0$ by assumption.  If $b$ is sufficiently positive for $M$, since $M_c=H^0(\PP,\widetilde{M}(c))$ for all $c\geq b$, it follows that $(H^1_{\mathfrak b}M)_c = 0$ for all $c\geq b$ as well.  In addition, $M_{\geq b}$ admits a linear resolution
\[
M_{\geq b}(b) \gets F_0\gets F_1\gets \dots \gets F_{p}\gets 0,
\]
where $\beta_{k,a}(M_{\geq b}(b))$ only if $a\geq 0$ and $|a|=k$.  We fix some $j=(j_1,\dots,j_t)\in \mathbb N^t$ with $|j|= i$.  To show that $b$ lies in the multigraded regularity, we need to show that $(H^{i+1}_{\mathfrak b} M)_{b-j}=0$ which amounts to showing that $H^i(\PP, \widetilde{M}(b-j))=0$. For this, it suffices to check that $H^{i+k}(\PP, \widetilde{F}_k(-j))=0$ for all $k$.   Since the resolution is linear, this amounts to checking that $H^{i+k}(\PP, \cO(-a-j))=0$ whenever $a,j\geq 0$ and $|a|=k$ and $|j|=i$, and this holds for line bundles on $\PP$.

Conversely, let $b$ lie in the multigraded regularity of $M$.  The first half of the second condition of Proposition~\ref{high truncation} is immediate from the definition of multigraded regularity~\cite[Definition~1.1]{maclagan-smith}.  It then suffices to show that $M_{\geq b+1^t-1_i}$ has a linear resolution, as the second half of the second condition will follow.  By Corollary~\ref{quadrant Tate}, it suffices to show that $\mathbf{R}(M_{\geq b+1^t-1_i})$ is acyclic.  Let $T$ be the Tate resolution of $\widetilde{M}$ and consider the corner complex $T_{\leq b-1_i}[-t]\to T_{\geq b+1^t-1_i}$.  Since the higher cohomology of $\widetilde{M}$ vanishes for all multidegrees $c\geq b-1_i$, it follows that the quadrant complex $T_{\geq b+1^t-1_i}$ equals $\mathbf{R}(M_{\geq b+1^t-1_i})$.  So we want to show that $T_{\leq b-1_i}[-t]$, has no terms in cohomological degrees $> |b-1_i|-t$, or equivalently that $T_{\leq b-1_i}$ has no terms in cohomological degrees $\geq |b|$.  

By Theorem~1.12, any such term would correspond to a nonzero cohomology group $H^p(\PP,\widetilde{M}(b-1_i -j))$ where $j\in \mathbb N^t$ and $p\geq |1_i+j|$.  However these groups are all zero since $b$ is in the multigraded regularity.
\end{proof}

\section{An example on $\PP^1 \times \PP^1$}\label{P1xP1 examples}

%\frank{This example before the Tate extension theorem is useful, because it explains a little the encoding of $T$ by its  cohomology table}
We recall that, by~\cite{cox}, we can present a coherent sheaf on $\PP^1\times \PP^1$ as a bigraded module over the Cox ring $S= K[x_{0}, x_{1}, y_{0}, y_{1}]$. Any sheaf on 
$\PP^{3}$ with a $G_{m}$ action will define such a module. For example, consider the universal sub-bundle $U$ on $\PP^3$ and the corresponding sheaf $\sF$ on $\PP^1\times \PP^1$.
As in~\cite[Theorem~4.1]{EFS}, the cohomology table of $U$ is
given by the Betti table of its Tate resolution:
$$
\begin{tabular}{r|rrrrrrrrrrrrr}
             & -5&  -4& -3& -2 & -1 & 0 & 1 & 2 & 3 & 4 & 5 \cr\hline
%      total:& 140& 84 &45& 20&  6& 1 \cr \hline
         3:  &  120&  70 & 36 & 15 & 4 &. &  .&  . & . & . & .\cr
         2: &  .&  . & . & . & . &. &  .&  . & . & . & .\cr
          1: &  .&  . & . & . & . &1 &  .&  . & . & . & . \cr
          0:&  .&  . & . & . & . &. & 6&20&45&84&140
\end{tabular} 
$$
Consider the map  $T^{0}(U)\to T^{1}(U)$ indicated by the numbers 1 and 6 in the table above.  This corresponds, in the $\ZZ^2$ grading, to the map 
$\omega_E \to \omega_E(-2,0)\oplus  \omega_E^4(-1,-1) \oplus \omega_E(0,-2)$
 defined by the matrix 
 $$
m=(e_0e_1,e_0f_0,e_1f_0,e_0f_1,e_1f_1,f_0f_1)^t
 $$
 where $V=V_1\oplus V_2=\langle e_0,e_1,f_0,f_1\rangle$.
By reciprocity, 
$$
\bL(\image\ m) \to M \to 0
$$
 is the minimal free resolution of the module  of global sections $M=\sum_{(a,b) \in \ZZ^2} H^0(\sF(a,b))$ of a  rank $3$ vector bundle $\sF$.

The cohomology table of $\sF$, written as a matrix over $\ZZ[h]$, is:
$$
\Big(\sum_{i=0}^2 {\rm h}^i(\sF(a,b))\cdot h^i \Big)_{-3 \le a,b \le 3}
%$$
%whose ranks are shown in the table
%$$
=\quad \begin{pmatrix}28 h&
{\color{blue}           18 h}&
{\color{red}      8 h}&
{\color{red}        2}&
      12&
      {\color{blue}22}&
      {\color{red}32}\\ % first line
{\color{red}       20 h}&
      13 h&
 {\color{blue}     6 h}&
{\color{blue}      1}&
{\color{red}        8}&
      15&
      {\color{blue}22}\\ % second line
 {\color{blue}      12 h}&
 {\color{red}      8 h}&
      4 h&
      0&
{\color{blue}       4}&
{\color{red}        8}&
      12\\ % thirdline
  {\color{black}      4 h}&
 {\color{blue}      3 h}&
 {\color{red}        2 h}&
      h&
      0&
{\color{blue}       1}&
{\color{red}       2}\\
{\color{black}       4 h^{2}}&
 {\color{blue}      2 h^{2}}&
   {\color{blue}    0}&
 {\color{red}       2 h}&
      4 h&
 {\color{blue}      6 h}&
{\color{red}        8 h}\\
 {\color{red}      12 h^{2}}&
{\color{black}       7 h^{2}}&
 {\color{blue}      2 h^{2}}&
  {\color{blue}     3 h}&
  {\color{red}      8 h}&
      13 h&
 {\color{blue}      18 h}\\
{\color{blue}       20 h^{2}}&
 {\color{red}      12 h^{2}}&
 {\color{black}      4 h^{2}}&
 {\color{black}      4 h}&
  {\color{blue}     12 h}&
 {\color{red}       20 h}&
      28 h\\
      \end{pmatrix}
      $$
Note that the corresponding Tate resolution 
$$
T=\bT(\sF): 
\cdots\to
%{\color{blue}T^{-5}}\to
{\color{red}T^{-4}}\to
{\color{black}T^{-3}}\to
{\color{blue}T^{-2}}\to
{\color{red}T^{-1}}\to
T^{0}\to
{\color{blue}T^{1}}\to
{\color{red}T^{2}}\to
{\color{black}T^{3}}\to
{\color{blue}T^{4}}\to
%{\color{red}T^{5}}\to\cdots
$$
has terms from the diagonal colored bands of the cohomology table  with $T^0$ corresponding to the main diagonal above. The maps have components
corresponding to arrows pointing to terms in north, north-west and west direction in the next colored band. 

For example, the  matrix $m$ above corresponds to the three arrows
encoded in the submatrix 
$$
\begin{pmatrix} {\color{blue} 1} && \\ & {\color{blue}4} & \\{\color{black} h} && {\color{blue}1}\end{pmatrix} .
$$
As another example, there is a $(3+2+2+3)\times (8+2+2+8)$ submatrix of ${\color{blue}T^{-2}}\to
{\color{red}T^{-1}}$
encoded by 
$$
\begin{pmatrix} {\color{red}8h} &&& \\ {\color{blue}3h}& {\color{red}2h} && \\ {\color{blue}2h^2}&& {\color{red}2h} &\\ &{\color{blue}2h^2}&{\color{blue}3h}& {\color{red}8h}\end{pmatrix}.
$$
Here $6$ of the $16$ blocks  are zero because of the north/west condition. 

In this example, the cohomology of $\sF$ is ``natural'', corresponding to the fact that each entry of the cohomology table is a monomial. In the more general case the terms of an entry that is not a monomial would contribute summands to different $T^d$.

The Tate resolution $T(U)$ on $\PP^3$ can also be thought of the
as the complex obtained by considering the
$\cornerT 0 (\sF)$ with respect to the natural ``coarse'' $\ZZ$ grading.

\begin{remark}\label{rmk:divergent examples}
In the above example, by considering the cohomology table of $\sF$, we can see that $(1,1)$ lies in the multigraded regularity of $M$.  However, since $H^1(\sF)=k^1$, this would yield a term of homological degree $1$ in the corner complex for $M$ at $(1,1)$. The total betti numbers of the corner complex $\cornerT {1^2} \sF$ in the coarse grading are
$$
\begin{tabular}{r|rrrrrrrrrrrrr}
&-3 &-2& -1 & 0 & 1 & 2 & 3 & 4  \cr\hline
%      total:& 140& 84 &45& 20&  6& 1 \cr \hline
         3:   & 36 & 15 & 4 &. &  .&  . & . & . \cr
         2:  & 10& 8& 6 &4 &  1&  . & . & . \cr
          1: &. & . & . &. &  .&  . & . & . \cr
          0: &. & . & . &. & 4&16&39&76
\end{tabular} 
$$
as one can read from the cohomology table of $\sF$ given above (using some values in addition to those shown in the table). In particular the complex 
$\bR(M_{\geq (1,1)})$, with betti numbers $4, 16, 39, 76,\dots$,  is not exact, as one sees from the ``1'' in the second row of the Betti table, which represents an element of the kernel  
of the differential 
$$
 \omega_E^{8}(-2,-1)\oplus \omega_E^{8}(-1,-2) \to \omega_E^{12}(-3,-1)\oplus \omega_E^{15}(-2,-2) \oplus \omega_E^{12}(-1,-3)
$$
 which is not in the image of $\omega_E^4(-1,-1)  \to \omega_E^{8}(-2,-1)\oplus \omega_E^{8}(-1,-2)$.
Thus $\bR(M_{\ge (1,1)})$ is not acyclic and $M_{\ge (1,1)}$ cannot have a linear resolution by Theorem~\ref{reciprocity}. Hence $(1,1)$ is not sufficiently positive for $M$ in our sense.  However, both $(1,2)$ and $(2,1)$ are sufficiently positive.
\end{remark}

\section{Injective and projective resolutions}\label{resolution matters}
To prepare for the proof of Theorem~\ref{Tate from BW}, we remind the reader of some general results about resolutions of complexes.

Let $M^{\bullet}:\cdots \to M^{i}\to M^{i+1}\to\cdots $ be a bounded-above complex of modules over a ring $R$.
A \emph{projective resolution} of $M^{\bullet}$ is  a complex of projective modules
$F^{\bullet}$ and a quasi-isomorphism $F^{\bullet} \to M^{\bullet}$, and similarly for
injective resolutions. Such resolutions were constructed in the famous book of Cartan-Eilenberg~\cite{Cartan-Eilenberg} by 
putting together
resolutions of the various kernels and cokernels of maps in $M^{\bullet}$, but the same
goal can be accomplished using iterated mapping cones. For the reader's convenience we give a proof
of this elementary result.

\begin{proposition}[Resolution of complexes by iterated mapping cones]\label{res of cplx}
Let $R$ be a ring. Let
$$
M^\bullet: \cdots \to M^{k-1} \to M^k \to \cdots
$$
be a bounded above complex of $R$-modules, and let
$$
M^{\bullet}_{p}: 0\to M^{p}\to M^{p+1}\to \cdots
$$ 
be the subcomplex of $M^{\bullet}$ obtained by truncation. For each $p$, let 
% $K^p\subseteq M^p$ be the kernel of $M^p\to M^{p+1}$ and let
$$
G_{p}^{\bullet}: \ldots \to G^{p-1}_{p} \to G^{p}_p \to M^p \to 0
$$
be a projective resolution of $M^p$.  
%$$
%G_{p}^{\bullet}: \ldots \to G^{p-1}_{p} \to G^{p}_p \to M^p \to 0
%$$
%be a projective resolution of $M^p$.  
There is a sequence of projective complexes 
$$
\cdots\subset F_{p+1}^{\bullet}\subset F_{p}^{\bullet}\subset \cdots
$$
 and surjective quasi-isomorphisms $F_{p}^{\bullet}\to M_{p}^{\bullet}$
such that $F^{\bullet}_{p}$ is the mapping cone of a map $G_{p}^{\bullet}[-1]\to F_{p+1}^{\bullet}$,
and $F^{\bullet} := \cup_{p} F_{p}^{\bullet}$ is a projective resolution of $M^{\bullet}$. Thus
the $k$-th term of $F^\bullet$ is $F^{k} = \sum_{p\leq k} F_{p}^{k}$. 

A similar result holds for bounded below complexes and injective resolutions.
\end{proposition}

If $\phi: A\to B$ is a map of complexes 
we write the mapping cone of $\phi$ as $[\phi]$ or $[A\rTo{\phi} B]$ or even $[A \to B]$ when $\phi$ is clear from context, so that there
is an exact sequence of complexes
$$
0 \to B \to [\phi] \to A[1] \to 0.
$$

\begin{proof}
Note that $M= \cup_{p} M^{\bullet}_{p}$. 
%Choose projective resolutions $G^{\bullet}_{i}$ of $M_{i}$.
%We will construct a sequence of complexes $F_{0}^{\bullet}\subset F_{1}^{\bullet}\subset\dots$ of
%projective modules of the form
%$$
%F^{\bullet}_{i} : \cdots \to F^{i-1}_{i}\to F^{i}_{i} \to 0
%$$ 
%and
%maps $F^{\bullet}_{i}\to M_{i}^{\bullet}$, such that $F^{i}_{i}\to M^{i}$ is surjective,
% that induce  isomorphisms on $j$-th cohomology modules for 
%all $j\geq i$. 
Since  $M^{\bullet}$ is bounded above there is a $k$ such that $M^{p} = 0$ for all $p>k$, and we may take $F_{p}^{\bullet} = M_{p}^{\bullet},$ the complex whose terms are all 0, for $p>k$. We now use descending induction on $p$.

Suppose that $F^{\bullet}_{p+1} \to M^{\bullet}_{p+1}$ is a projective resolution. 
We will show that the map of complexes
$$
\phi: \bigl(0\to M^{p}\to 0\bigr)[1] \longrightarrow\bigl(0\to M^{p+1}\to M^{p+2}\to\cdots  \bigr)
$$
induced by the differential of $M^{\bullet}$  lifts to a map
 $\phi': G_{p}^{\bullet}[-1] \to F_{p+1}^{\bullet}$, so that the maps $\phi$ and $\phi'$ are quasi-isomorphic,
 and we define $F_{p}^{\bullet}$ to be the mapping cone $F_{p}^{\bullet} := [\phi']$.
% We have the short exact sequence for the mapping cone $0\to F_{p+1}^\bullet \to F_p^\bullet \to G^p_\bullet \to 0$.
From the long exact sequence of the mapping cone, we obtain
$H^i F_{p}^{\bullet}\cong H^i F_{p+1}^\bullet\cong H^i M^\bullet$ for $i>p+1$, as well as a four-term exact sequence
\[
0\to H^p F_p^\bullet \to H^p G_p^\bullet \to H^{p+1} F_{p+1}^\bullet \to H^{p+1}F_{p}^\bullet\to 0.
\]
Since $H^p G_p^\bullet = M^p$ and $H^{p+1} F_{p+1}^\bullet= \ker(M^{p+1}\to M^{p+2})$ we immediately obtain that $H^i F_p^\bullet = H^i M^\bullet_p$ for $i=p,p+1$.
Hence
$F_p^\bullet$ is
a projective resolution of $M^{\bullet}_p$.
By construction, $F_{p+1}^{\bullet}\subset F_{p}^{\bullet}$ and they agree in  degree $\geq p+1$.
It follows that $F^{\bullet}$ is a projective resolution of $M^{\bullet}$, as claimed.
 
It remains to produce the map $\phi'$.
The module $Z$ of $(p+1)$-cycles in $M_{p+1}^{\bullet}$  is  the same as that in $M_{p}^{\bullet}$, and 
$M_{p+1}^{\bullet}$ has no $(p+1)$-boundaries, so the module $\tilde Z$ of $(p+1)$-cycles
in $F_{p+1}^{\bullet}$ maps surjectively to $Z$. Since $M^{q}_{p+1} = 0$ for $q< p+1$, the quotient complex
of $F^{\bullet}_{p+1}$ obtained factoring out $F_{p+1}^{r}$ for $r\geq p+1$ is a resolution of the kernel of
$\tilde Z \to Z$.
Since $G^{\bullet}_{p}$ is a projective resolution of $M^{p}$, the
map $G_{p}^{p} \to M^{p}$ lifts to a map  $G_{p}^{p}\to \tilde Z$, and we continue lifting to obtain
the desired map of complexes $\phi':G_{p}^{\bullet}[-1]\to F_{p+1}^{\bullet}$, as in the diagram:
$$
\begin{diagram}[small]
&&& && F^{p+2}_{p+1}&\rTo&M^{p+2}\\
&&& && \uTo&&\uTo\\
&&& && F^{p+1}_{p+1}&\rTo&M^{p+1}\\
&&& &\ruTo& \cup &&\cup \\
&F^{p-1}_{p+1}&\rTo&F^{p}_{p+1} &\rTo& \tilde Z &\rTo &Z &\rTo & 0\\
\phi':&\uTo &&\uTo&&\uTo&&\uTo\\
&G_{p}^{p-2}&\rTo&G_{p}^{p-1}&\rTo&G_{p}^{p}&\rTo& M^{p}&\rTo &0.
\end{diagram}
$$
\end{proof}

%\begin{proposition}[Projective/injective resolution of complexes]\label{res of cplx}
%Let $R$ be a ring and
%$$
%M^\bullet: \ldots \to M^{k-1} \to M^k \to 0 
%$$
%be a bounded above complex of $R$-modules. For each $p$, let 
%$$
%F^p_\bullet: \ldots \to F^p_1 \to F^p_0 \to M^p \to 0
%$$
%be a projective resolution of $M^p$.  An iterated mapping cone between the $F^{p}_\bullet$ is a projective resolution 
%$$
%F^\bullet: \ldots \to F^{k-1} \to F^k \to 0
%$$
%of $M^\bullet$ with terms
%$$
%F^\ell =\bigoplus_{p-q=\ell \atop p \le k, q\ge 0} F^p_q.
%$$
%A similar result holds for bounded below complexes and injective resolutions.
%\end{proposition}

The reciprocity theorem for resolutions of modules over $E$ and $S$ is a special case of a reciprocity theorem
for complexes, proved in the same way.

\begin{theorem}\label{Reciprocity for complexes} Let $M^\bullet$ be a bounded complex of finitely generated $S$-modules and $P^\bullet$ a bounded complex of finitely generated $E$-modules.
Then $\bR(M^\bullet)$ is an injective resolution of the complex $P^\bullet$ if and only if $\bL(P^\bullet)$ is a projective resolution of $M^\bullet$.
\end{theorem}

\begin{proof} The key point is that $\bL$ and $\bR$ are adjoint functors. See \cite[Theorem 2.6]{EFS}, \cite[Theorem 3]{BGG} for details.
\end{proof}

\section{Tate Resolutions from Beilinson Monads}\label{Tate from Monad}

Let $\sF$ be a finite complex of sheaves.
Generalizing the case of a single sheaf,
we may construct a \emph{Tate resolution $\bT(\sF)$} as follows.
We may represent $\sF$ by a finite complex of graded $S$-modules $M^\bullet$.
Since $\bR$ is a functor, we may apply it term by term to the
complex $M^\bullet$ to get a double complex of modules over $E$. We then
can apply the procedure of Section~\ref{construction}. We write $M^{\bullet}_{\geq c}$
for the complex
$$
\cdots \to M^{p}_{\geq c}\to M^{p+1}_{\geq c}\to \cdots,
$$
where $M^{p}_{\geq c}$ is the truncation in degrees $\geq c$ as in Section~\ref{construction}.

\begin{theorem}[Tate resolution and Beilinson representative of a complex]\label{Tate from BW}
Let $\sF$ be a bounded complex of coherent sheaves. Then there exists a  minimal Tate resolution $T=\bT(\sF)$ such that
$\bU(T(c)[c])$ is quasi-isomorphic to  $\sF(c)$ for every $c \in \ZZ^t$.  In particular $\bT(\sF)$ satisfies Equation~\ref{eqn:hypercohomology}.
\end{theorem}

\begin{proof}
 Let $\sF$ be represented by a finite complex $M^\bullet$ of finitely generated graded $S$-modules, and let $c\in \ZZ^{t}$ be sufficiently positive for each
of the finitely many nonzero modules $M^{p}$ in $M^{\bullet}$. Each module 
$M^{p}_{\geq c}$ has a linear minimal free resolution $L^{p}$ and the differentials of $M^{\bullet}$
lift to homogeneous maps $L^{k}\to L^{k+1}$. Since adding a homotopy would make the maps inhomogeneous, the lifted maps are unique and compose to 0. We may write
$L^{k} = L(P^{k})$ for suitable finitely generated $E$-modules $P^{k}$, and
it follows that the lifted maps on the $L^{k}$ make the $P^{k}$ into a complex, which we denote $P^{\bullet}$.

Since $\bR$ is a 
functor, we may regard $\bR(M^{\bullet})$ as a double complex.
By Theorem~\ref{reciprocity} the complex $\bR(M^k_{\ge c})$ is the injective
resolution of 
$
P^k.
$
We regard the differentials in $\bR(M^{\bullet})$
coming from the individual $\bR(M^{k}_{\geq c})$ as the ``horizontal differentials''
in this double complex. With this convention the ``horizontal homology''
of $\bR(M^{\bullet}_{\geq k})$ is the complex $P^\bullet$.

Replacing $c$ by $c+1^{t}$ if necessary, we may ensure that $c-1^{t}$ is sufficiently
positive for each $M^{k}$. Then the minimal projective resolution of $P^{k}$
starts with the term $\bR(M^k_{c-1^t})$ by Theorem \ref{cut complex1}.
%, and the composite map
%$$
%\bR(M^k_{c-1^t}) \to P^{\bullet} \to \bR(M^{\bullet}_{\geq k}
%$$

Let $\Tail_c (M^\bullet) \to P^\bullet$ be a minimal projective resolution. By 
Theorem~\ref{res of cplx}  we can obtain such resolution by minimizing an iterated mapping cone of projective resolutions of the $P^k$. As in the case of a single module, we set
$$
\bT(\sF) := \left({\lim_{\infty\leftarrow c}}' \Tail_{c}(M^\bullet)\right)[1-t].
$$
Since $\bL(P^\bullet) \to M^\bullet_{\ge c}$ is a  resolution, the result follows as in Section \ref{monads}. 
\end{proof}

Conversely given any locally finite complex $T$ of free $E$-modules, we get a finite complex of sheaves
$U=\bU(T)$, which we may regard as a Beilinson representative of an 
object $\sF$ in $D^b(\PP)$. We will now show how to construct the Tate resolution
$\bT(\bU(T))$ directly, in the context of resolutions over the exterior algebra,
 without passing through the category of sheaves. We will achieve this goal in several steps.

There is a unique smallest free subquotient complex $\bW(T)$ of $T$,
called the \emph{Beilinson window} of $T$, such that $\bU(T)=\bU(\bW(T))$. Indeed, if 
$$
T^d = \sum_{a \in \ZZ^t} B_{a} \otimes \omega_E(-a)
$$ 
$\bW(T)$ is the subquotient complex with terms
$$
BW^d= \sum_{a \in \ZZ^t \atop 0 \ge a \ge -n} B_{a} \otimes \omega_E(-a).
$$

Now suppose that $BW$ is any complex of finitely generated free $E$-modules 
that are direct sums of modules of the form $\omega_{E}(b)$ with $0 \leq b \leq n$.
We call $\bU(BW)$ \emph{minimal} if $BW$ is minimal. In general,
$\bU(BW)$ is the direct sum of a minimal complex $U_{min}$ and
trivial complexes of type $0 \to U^b \to U^b \to 0$ for various
degrees $b$ with $0 \le b \le n$ and various shifts; this follows
from the corresponding fact for $BW$. 
  
%\begin{theorem}\label{Tate from Beilinson}
%Let $BW$ be a finite minimal complex of the form $BW^{d} = \sum_{0 \ge a \ge -n} B_{a}^{d}\otimes \omega_E(-a)$
%and let
%$P(BW)$ be the subcomplex
%$$
%P(BW) =\langle f \in BW \hbox{ homogeneous } \mid \deg f \not\ge 1^t \rangle \subset BW.
%$$
%Let $BW/P(BW) \to T^+$ be a minimal 
%injective resolution,  let $T^- \to P(BW)$ be a minimal projective resolution, and set
%$$
%cT:=[\, [T^-[-1] \to BW] \to T^+[1]\,].
%$$
%The complex $cT$ is isomorphic to the corner complex $(\cornerT {1^t})[-1]$ of a Tate resolution $T$ with $BW=\bW(T)$. 
%$\end{theorem}

\begin{theorem}\label{corner complex from Beilinson} Let $T$ be a minimal Tate resolution and $BW=\bW(T)$ its Beilinson window.
Let
$P(BW) \subset BW$ be the subcomplex
$$
P(BW) =\langle f \in BW \hbox{ homogeneous } \mid \deg f \not\ge 1^t \rangle \subset BW.
$$
If $BW/P(BW) \to I$ is a minimal 
injective resolution and  $F \to P(BW)$ is a minimal projective resolution, then the corner
complex $\cornerT {1^t}$ is isomorphic to
$$
[\, [F \to BW] \to I\,].
$$
and the shifted corner complex $cT :=  (\cornerT {1^t})[-1]$ has $BW$ as a subquotient complex.
\end{theorem}

\begin{example}\label{corner of terms}
 If $BW$ consist of a single term $\omega_E(b)$ with $0 \le b \le n$ then $\bU(BW) = U^b$ and the corner complex
$(\cornerT {1^t})[-1]$ as described in Example \ref{special corner complexes} for $T=\bT(U^b)$ coincides
with $cT$. Indeed, in this case
\begin{align*}
P(BW) &:=  \langle f \in \omega_E(b) \hbox{ homogeneous }\mid \deg f \not\ge 1^t \rangle\\
         &=  \langle f \in \omega_E(b) \hbox{ homogeneous }\mid \exists j : (\deg f)_j \le 0 \rangle \\
         &= \ker (\omega_E(b) \to \tensor_{j=1}^t \Lambda^{b_j+1} W_j \tensor \omega_E(b)).
\end{align*}
Here we consider the vector space $W_{j}$ to be concentrated in degree $-1_j$, so 
$\Lambda^{b_j+1} W_j \tensor \omega_E(b))$ is isomorphic to a direct sum of copies of $\omega_{E}(-1^t)$.
\end{example}

\begin{proof}[Proof of  Theorem~\ref{corner complex from Beilinson}.]
The lower quadrant $qT=T_{\leq 0}$ is a quotient complex of $cT$ and the shifted upper quadrant $uT=T_{\ge 1^t}[t-1]$ is a subcomplex of $cT$, which in turn is a cone over the corner map between these two quadrant complexes. $BW$  is a subcomplex of $qT$, and we let
$$
qT^{-}  :=qT /BW \cong  \sum_{a \in \ZZ^t\atop 0 \le -a, -a \not\leq n} B_{a} \otimes \omega_E(-a) 
$$
denote the part of $qT$ outside the Beilinson window. 
Since
$$
\Lambda^{n_i+2} V_i =0,
$$
the corner map induces the zero map from $qT^{-}$ to  $uT$.

Thus
$$
cT[1] = [\, [ qT^{-}[-1] \to BW] \to uT \,]=[ qT^{-}[-1] \to [BW \to uT\,]\,]
$$
where $[A \to B]$ denotes the cone over a map of complexes as in Section~\ref{resolution matters}. 

By the exactness of $cT$ we have 
$$
P(BW) =\langle f \in BW \hbox{ homogeneous } \mid \deg f \not\ge 1^t \rangle = \ker (BW \to uT).
$$
The last equality holds since, by Example \ref{corner of terms}, it holds for every term $\omega_E(b)$ of $BW$, and because the injective resolution can be obtained  by minimizing an iterated mapping cone as in Proposition~\ref{res of cplx}.
Since $T$ and hence $cT$ are minimal, we recover $uT$ from $BW$ as the minimal injective resolution of $BW/P(BW)$. Also $F=qT^-[-1]$ is the minimal projective resolution of 
$P(BW) \subset BW$. Hence we recover $qT=[F \to BW]$ and $cT[1]=[qT \to uT]$ as cones.
\end{proof}

\begin{algorithm}\label{lowerQuadTateExtension} [lower quadrant of the Tate extension]

\noindent
{\bf Input:} $BW\colon 0 \to BW^r \to \ldots \to BW^s \to 0$, a finite minimal complex with terms
$BW^{d} = \sum_{0 \ge a \ge -n} B_{a}^{d}\otimes \omega_E(-a)$ and 
 an integer $d$ with $d < r$.\\
{\bf Output:} The lower quadrant complex $$qT=T_{1^t}(\{1,\ldots t\}, \emptyset,\emptyset)=[F\to BW]$$ of Theorem \ref{corner complex from Beilinson}
in homological degrees $d, \ldots, s$.
\begin{enumerate}
\item Set $k:=s, \ qT^k :=BW^s$ and $qT^{k+1}:=0$.
\item While $k > d$ do:
\begin{enumerate}
\item Compute $L^k=\ker(qT^k \to qT^{k+1})/\image(BW^{k-1} \to qT^k)$.
\item Compute minimal generators of 
$$P(L^k)= \langle f \in L^k \hbox{ homogeneous}) \mid \deg f \not\ge 1^t \rangle$$. 
\item Choose a map
$ F^k \to \ker (qT^k \to qT^{k+1})$ from a free $E$-module of the form $F^k = \sum_{0 \ge a \not\ge -n} B^k_a \tensor \omega_E(-a)$
 such that the  image of the generators in $F^k$ generates $P(L^k)$ minimally.
\item Set $qT^{k-1} := BW^{k-1} \oplus F^k$, and extend the differential of $qT$ to include the map $qT^{k-1} \to qT^k$.
\item Replace $k$ by $k-1$. 
\end{enumerate}
\item Return the complex 
$$ qT^d \to qT^{d+1} \ldots \to qT^s \to 0.$$
\end{enumerate}
\end{algorithm}

\begin{proof} We keep the notation of the proof of Theorem \ref{corner complex from Beilinson}. Representatives in  $qT^k$ of elements of $P(L^k)$ map to zero in $uT[1]$ for degree reasons. By the exactness of $cT$ they must be covered by elements in $qT^-$. The result follows.
\end{proof}

\begin{algorithm}\label{corner complex computation}[Corner Complex]

{\bf Input:} $BW\colon 0 \to BW^r \to \ldots \to BW^s \to 0$, a finite minimal complex with terms
$BW^{d} = \sum_{0 \ge a \ge -n} B_{a}^{d}\otimes \omega_E(-a)$ and 
 a degree $a \in \ZZ^t$  with $a <-n$.
{\bf Output:} A finite piece of the corner complex $\cornerT a$ of the Tate resolution $T=\bT(\bU(BW))$.
\begin{enumerate}
\item Set $d= r +|a|-|n|-t$ and $e=r+|a|$
\item Apply Algorithm \ref{lowerQuadTateExtension} to compute $ qT: qT^d \to qT^{d+1} \ldots \to qT^s \to 0.$
\item Collect $\ell T$, $uT$ the terms and maps of $qT$ contributing to the lower part $T_a(\{1,\ldots,t\},\emptyset,\emptyset)$ and the upper part $T_a(\emptyset,\emptyset,\{1,\ldots,t\})$ respectively.
\item Compute the corner maps $\ell T[t-1] \to uT$ and form the subquotient complex $cT=[\ell T[t-1] \to uT]$ of $\cornerT a$.
\item Then the differential $cT^{e-1} \to cT^{e}$ and $(\cornerT a)^{e-1} \to (\cornerT a)^{e}$ coincide. 
\item Extend the part of the corner complex to the complex
$$(\cornerT a )^{e-1} \to \ldots \to (\cornerT a)^{s}$$
by injective (and possibly projective) resolutions of the single correct differential from Step 5, and return the result.
\end{enumerate}
\end{algorithm}

\begin{proof} Since the nonzero summands $B^r_c \tensor \omega_E(-c)$ of $BW^r$ satisfy $c \le 0$ taking syzygy we see that the nonzero summands $B^{e-1}_c \tensor \omega_E(-c)$ of $qT^{e-1}$ satisfy $c \not\ge  a$. Thus $qT^e$ is the last term which possibly could contribute to $uT$ and $qT^d$ is the last term which possibly maps to $uT^e$ via a corner map. Moreover since $a<-n$ the lower quadrant $T_a(\{1,\ldots,t\},\emptyset,\emptyset)$ maps to no term $B^\ell_c \tensor \omega_E(c)$ of $T^\ell$ with $c \not\le 0$. Thus $cT^{e-1} \to cT^{e}$ and $(\cornerT a)^{e-1} \to (\cornerT a)^{e}$ coincide. 
\end{proof}

\begin{algorithm}\label{cohomologyTable} [Cohomology Table]

\noindent
{\bf Input:}  $BW\colon 0 \to BW^r \to \ldots \to BW^s \to 0$, a finite minimal complex with terms
$BW^{d} = \sum_{0 \ge a \ge -n} B_{a}^{d}\otimes \omega_E(-a)$; and
two  degrees $a, b \in \ZZ^t$  with $a <-n$ and $0 < b$.

\noindent
 {\bf Output:} The cohomology table of $\sF =\bU(BW)$
 $$ \{ c \in \ZZ^t \mid a \le c \le b \} \to \ZZ[h,h^{-1}]$$
 $$ c \mapsto \sum_{k\in \ZZ} \dim \HHH^k( \sF(c)) h^k$$ 
 %\begin{enumerate}
 %\item Use Algorithm~\ref{corner complex computation} to compute the differential 
 %$$(\cornerT a)^{d+|n|} \to (\cornerT a)^{d+|n|+1}$$
 %of the corner complex $\cornerT a$ of $T=\bT(\sF)$.
 %\item Set $k={d+|n|}$ and $ \partial_k: uT^k \to uT^{k+1}=(\cornerT a)^{k} \to (\cornerT a)^{k+1}$.
 %\item {\bf repeat} 
 %\begin{enumerate}
 %\item replace $k$ by $k+1$
 %\item Compute a minimal free injective hull $uT^{k+1}:= \sum_{e\in \ZZ^t} B^{k+1}_e \tensor \omega_E(-e)$ of $\coker \partial_{k-1}$ and define $ \partial_{k}\colon uT^{k} \to uT^{k+1}$ as the induced map.
% \item Set $$uT^{k+1} := \sum_{e\in \ZZ^t \atop e \le b} B^{k+1}_e \tensor \omega_E(-e)$$ and define $ \partial_{k}\colon uT^{k} \to uT^{k+1}$ as the induced map.
% \end{enumerate}
 %{\bf until} $\sum_{e\in \ZZ^t \atop e \le b} B^{k+1}_e \tensor \omega_E(-e)=0$.
 %\item Read off the cohomology table from the complex
 %$$ uT^{d+|n|} \to \ldots \to uT^k$$
 %by collecting the dimensions 
 %$$
 %\sum_{\ell=d+|n|}^k \dim B^\ell_c \, h^{\ell-|c|} \in \ZZ[h,h^{-1}]
 %$$ 
 %for all $c \in \ZZ^t$ with $a \le c \le b$. 
 %\end{enumerate}
 \begin{enumerate}
 \item Set $e=r+|a|$
 \item Use Algorithm~\ref{corner complex computation} to compute the differential 
 $$(\cornerT a)^{e-1} \to (\cornerT a)^{e}$$
of the corner complex $\cornerT a$ of $T=\bT(\sF)$.
 \item Set $k={e-1}$ and $ \partial_k: uT^k \to uT^{k+1}=(\cornerT a)^{k} \to (\cornerT a)^{k+1}$.
 \item {\bf repeat} 
 \begin{enumerate}
 \item replace $k$ by $k+1$
 \item Compute a minimal free injective hull $I=\sum_{c\in \ZZ^t} B^{k+1}_c \tensor \omega_E(-c)$ of $\coker \partial_{k-1}$.
 \item Set $$uT^{k+1} := \sum_{c\in \ZZ^t \atop c \le b} B^{k+1}_c \tensor \omega_E(-c)$$ and define $ \partial_{k}: uT^{k} \to uT^{k+1}$ as the induced map.
 \end{enumerate}
 {\bf until} $uT^{k+1}=0$.
 \item Read off the cohomology table from the complex
 $$ uT^{e} \to \ldots \to uT^k$$
 by collecting the dimensions 
 $$
 \sum_{\ell=e}^k \dim B^\ell_c \, h^{\ell-|c|} \in \ZZ[h,h^{-1}]
 $$ 
 for all $c \in \ZZ^t$ with $a \le c \le b$. 
 \end{enumerate}
\end{algorithm}

\begin{proof} Let $T=\bT(\bU(BW))$.  The upper quadrant $uT =T_a(\emptyset,\emptyset, \{1,\ldots,t\})$ is a subcomplex of the corner complex, which by exactnessss can be computed from any of its differentials. Since the nonzero summands $B^r_c \tensor \omega_E(-c)$ of $BW^r$ satisfy $c \le 0$ taking syzygy we see that the nonzero summands $B^{e-1}_c \tensor \omega_E(-c)$ of $uT^{e-1}$ satisfy $c \not\ge  a$. So these do not contribute to the desired cohomology table. For the computation in the injective direction we may drop terms $B^\ell_c \tensor \omega_E(c)$ with $c \not\le b$, because the differential of the $T$ and the corner complex
goes in the cohomology table in northeast direction, so that these terms cannot contribute to terms or computations of terms in the desired cohomological range.
Hence the correctness follows from Equation~\ref{eqn:hypercohomology} applied to the complex $\bU(BW)$.
\end{proof}

\begin{corollary}\label{reconstruction} Starting from a Beilinson window $BW=\bW(T)$ of a Tate resolution,  one can recover arbitrarily large finite parts of the Tate resolution $T$ with computations over the exterior algebra.
\end{corollary}

\begin{proof} Apply the three algorithms above. 
\end{proof}

Using this result we can answer the question: ``When is a complex of sheaves on $\PP$ quasi-isomorphic to a vector bundle in homological degree 0?''

\begin{corollary}\label{detect sheaves}
Let $T$ be a Tate resolution. Then $\bU(T)$ is a monad of a sheaf if and only if the cohomology table of a positive quadrant $T_c(\emptyset, \emptyset, \{1,\ldots, t \})$ for some $c \gg 0$ has only $h^0$ entries. It is a monad of a vector bundle if only if in addition also the cohomology table of every very negative quadrant $T_c(\{1,\ldots, t \},\emptyset,\emptyset )$ for $c \ll 0$ has only $h^{|n|}$ entries. \hfill \qed
\end{corollary}

The corollary above is not an effective criterion, because we have no algorithm to compute the whole Tate resolution. However as in the case of a single projective space, see for example \cite[Example 7.3]{EFS}, there is  an effective method for proving that a bounded complex is quasi-isomorphic to a vector bundle.

\begin{proposition} \label{detect vector bundles} Consider $a,b\in \ZZ^t$ with $b>0$ and let $t_0<t_1 \in \ZZ$. Let $\sF \in D^b(\PP)$. Suppose
$$\HHH^{*} \sF(a+bt) = \HHH^{|n|} \sF(a+bt) \hbox{ for }  t \in [t_0,t_0+|n| ]$$ and  $$\HHH^* \sF(a+bt) = \HHH^0 \sF(a+bt) \hbox{ for } t \in [t_1-|n|, t_1]. $$
Then $\sF$ is quasi isomorphic to a vector bundle on $\PP$.
\end{proposition}

\begin{proof} Consider $\sL= \sO(b_1,\ldots, b_t)$. The line bundle $\sL$ is very ample and defines an embedding $\iota: \PP \hookrightarrow \PP^N$
with $N+1= \prod_{j=1}^t {b_j+n_j \choose n_j}$.  The cohomology table of $\sG=\iota_* \sF(a)$ with respect to $\sO_{\PP^N}(1)$ can be read of from the values of the cohomology table of $\sF$ along the 
integral line $ \ZZ \to \ZZ^t, t \mapsto a+bt$. If $\pi: \PP \to \PP^{|n|}$ denotes a linear Noether normalization of $\iota(\PP) \subset \PP^N$, then $\sG$ and $\pi_* \sG$
have the same cohomology table, and $\sG$ is a vector bundle, i.e. Cohen-Macaulay sheaf of dimension $|n|$, iff and only if $\pi_* \sG$ is a vector bundle.
The assumption $\HHH^* \sF(a+bt) = \HHH^0 \sF(a+bt) \hbox{ for } t \in [t_1-|n|, t_1]$ implies that for $\bT(\pi_* \sG)$ the assumption of \cite[Lemma 7.4]{EFS} is satisfied at position $t_1$. So in particular $\HHH^k_* (\pi_* \sG)= 0$ for $k >|n|$.
Similarly using the dual complex, $\HHH^{*} \sF(a+bt) = \HHH^{|n|} \sF(a+bt) \hbox{ for }  t \in [t_0,t_0+|n| ]$ implies using the Lemma, $\HHH^k_* (\pi_* \sG)= 0$ for $k < 0$.
Moreover, the two applications of the Lemma imply that the intermediate cohomology groups $\HHH^k_*(\sG)= \oplus_{d \in \ZZ} \HHH^k(\pi_* \sG(d))$ for 
$0< k < |n|$ have finite length. 
Thus $\pi_*\sG$  on $\PP^{|n|}$, and hence $\sF$ on $\PP$, is quasi-isomorphic to a vector bundle
by the Auslander-Buchsbaum formula~\cite{Eisenbud1995}.\end{proof}

\begin{remark} We have implemented these algorithms in our Macaluay2 package
\href{http://www.math.uni-sb.de/ag-schreyer/joomla/index.php/computeralgebra}{TateOnProducts.m2}~\cite{TateOnProducts}.
\end{remark}

%%%%%%%%%%%%%%%%%%%%
%%%%%%%%%%%%%%%%%%%%
\section{Application to split vector bundles}\label{splittings}
%%%%%%%%%%%%%%%%%%%%
%%%%%%%%%%%%%%%%%%%%
%%%%%%%%%%%%%%%%%%%%

Horrocks' splitting criterion (\cite{horrocks, OSS}) says that if a vector bundle on $\PP^n$ has no intermediate cohomology, then that vector bundle splits as a sum of line bundles.
%We consider this statement from a slightly different perspective.  First, we note that the cohomology of every vector bundle can be written as a sum of the cohomology tables of supernatural bundle~\cite[Theorem~0.5]{ES2008}.  Second, we note that the only supernatural bundles without intermediate cohomology are sums of line bundles of the form $\cO(-c)^{m}$.  Hence, Horrocks' criterion can be understood as the following principle: 
That is: {\em if the cohomology table of a vector bundle looks like a sum of line bundles, then the vector bundle itself splits as a sum of line bundles.}  

\begin{question}
Suppose that the cohomology table of a vector bundle $\cE$ on $\PP:=\PP^{n_1}\times \PP^{n_2}\times \dots \times \PP^{n_t}$ can be written as a positive integral sum of the cohomology tables of line bundles on $\PP$. Is $\cE$  a direct sum of line bundles?
\end{question}

We prove that the answer is ``yes'' under an additional hypothesis that is automatically satisfied in the case $t=1$.
For $c\in \ZZ^s$ we write $\gamma_{i,c}(\cE):= h^i(\PP, \cE(c))$ and we write $\gamma(\cE)=(\gamma_{i,c}(\cE))_{i,c}$ for the cohomology table of $\cE$. As throughout this paper, we use the termwise partial order for comparing integer vectors in $\ZZ^t$.

\begin{theorem}\label{thm: splitting theorem}
Let $\cE$ be a vector bundle on $\PP:=\PP^{n_1}\times \PP^{n_2}\times \dots \times \PP^{n_t}$ such that the cohomology table of $\cE$ decomposes as a positive sum of line bundles:
\[
\gamma(\cE)=\sum_{i=1}^s \gamma(\cO(c^{(i)})^{m_i}) \text{ where } c^{(i)}\in \ZZ^t.
\]
If $c^{(1)}\geq c^{(2)}\geq \dots \geq c^{(s)}$ then $\cE$ splits as
\[
\cE \cong \bigoplus_{i=1}^s \cO(c^{(i)})^{ m_i}.
\]
\end{theorem}
The proof uses our construction of the Beilinson monad.  We begin with a more general lemma.  
%We let $W=\{a \in \ZZ^t | (0,\dots,0)\geq \alpha \geq (-n_1,\dots, -n_t) \} \subseteq \ZZ^t$ be the set of vectors in $\ZZ^t$.
% such that the functor $U^{-a} \ne 0$.
%
%\david{the bracket notation in the following conflicted with our usage of $[**]$ for
%the mapping cone, and with usage of complexes elsewhere, so I changed it. }

\begin{lemma} \label{lem: splitting lemma}
Suppose that the cohomology table of $\cE$ splits as $\gamma(\cE)=\gamma(\cO^{m_1})+ \gamma(\cE')$, where $\gamma(\cE')$ is the cohomology table of some vector bundle $\cE'$, and where $\bU^{-1}(\cE')=0$.
Then $\cE$ splits as
\[
\cE = \cO^{m_1}\oplus \cE'',
\]
for some vector bundle $\cE''$ satisfying $\gamma(\cE'')=\gamma(\cE')$.
\end{lemma}
\begin{proof}
By assumption, $\bU^{0}(\cE)$ has $\bU^{0}(\cO^{m_1})=\cO^{m_1}$ as a summand,
so we may write $\bU^{0}(\cE)=\cO^{m_1}\oplus  \bU^{0}(\cE')$.  Note that $\bU^{-1}(\cE)
=\bU^{-1}(\cO^{m_1})\oplus \bU^{-1}(\cE')=0$. Thus $\bU(\cE) $ has the form
\[
= \cdots \to 0\to \cO^{m_1}\oplus \bU^0(\cE') \to \bU^1(\cE) \to \dots 
\]
Further, the summand $\cO^{m_1}$ must map to 0 in $\bU^1(\cE)$: minimality ensures that any map between sums of $\cO$ are zero, and the fact that the terms in the image of the functor $\bU$ form a strong exceptional collection implies that $\cO = U^{0}$ admits no nonzero maps to any $U^{a}$ with $a\neq 0$.  Thus, $\cO^{m_1}$ is a direct summand of the zeroth cohomology $\cE = H^{0}(\bU(\cE))$ as claimed.
%\[
%\Omega(\cE)=\cO^{m_1}\oplus \Omega(\cE')
%\]
%and hence $\cE=\cO^{m_1}\oplus \Omega(\cE')$.
\end{proof}

\begin{proof}[Proof of Theorem~\ref{thm: splitting theorem}] By induction, it suffices
to show that $\cO(c_{1})^{m_{1}}$ is a summand of $\cE$. 
Without loss of generality, we may assume that $c^{(1)} >c^{(2)}$.  Replacing $\cE$ by $\cE(-c^{(1)})$, we can further assume that $c^{(1)}=0$ and thus that $c^{(i)}<0$ for all $i=2,\dots,s$.  We will complete the proof by verifying the hypotheses of Lemma~\ref{lem: splitting lemma}.  

%\david{The following is hard to follow. If we're going to use the Lemma with a statement in terms of $\bT$, why not state it that way to begin with --- or at least say explicitly what the equivalent form is.}

Let $\cE':=\bigoplus_{i=2}^s \cO(c^{(i)})^{ m_i}$.  Then we have $\gamma(\cE)=\gamma(\cO^{m_1})+ \gamma(\cE')$ by assumption, and we need to show that $\bU^{-1}(\cE')=0$.
Since
\[
\bU^{-1}(\cE') = \bigoplus_{i=2}^s \bU^{-1}(\cO(c^{(i)})^{ m_i}),
\]
it suffices to show that $\bU^{-1}(\cO(b)) = 0$ for any $b<0$.  By Theorem~\ref{Terms of the Beilinson monad} we have
\[
\bU^{-1}(\cO(b)) = \bigoplus_{p+|a|=-1} U^{-a}\otimes H^p(\PP,\cO(b+a))
\]
Assume that $H^p(\PP,\cO(b+a))\ne 0$ for some $-n \leq a \leq 0$.  Then there exist $j_1,\dots,j_t$ such that $\sum_i j_i=p$ and
\[
H^p(\PP,\cO(b+a)) = \bigotimes_{i=1}^t H^{j_i}(\PP^{n_i},\cO(b_i+a_i)).
\]
There are two possibilities for each $j_i$:  either $j_i=0$ or $j_i=n_i$.  If $j_i=0$, then since $b_i$ and $a_i$ are both nonpositive, we must have $b_i=a_i=0$ and $j_i+a_i=0$.  If $j_i=n_i$, then since $a_i\geq -n_i$, we must have $j_i+a_i \geq 0$.  This implies that $p+|a|=\sum_i j_i+a_i \geq 0$ whenever $U^{-a}$ and $H^p(\PP,\cO(b+a))$ are both nonzero, and it follows that $\bU^{-1}(\cO(b))=0$ as claimed.
We may now apply Lemma~\ref{lem: splitting lemma} to split off a copy of $\cO^{m_1}$, completing the proof.
\end{proof}

\section{More Open Questions}

\begin{question}
Consider the case $t=2$. If $\sF \cong \sF_1 \boxtimes \sF_2$ then its Tate resolution has  the structure of a double complex:
$$ \bT(\sF) \cong \bT(\sF_1) \tensor_K \bT(\sF_2) \hbox{ as complex of } E = E_1 \tensor_K E_2 \hbox{ modules}.$$
Similarly a direct sum of box-products has a Tate resolution which is naturally a double complex.
Does the converse hold? 
In other words, assume that $\bT(\sF)$ admits the structure of a double complex, so that the differential $\partial$ decomposes as $\partial=\partial_h+\partial_v$, and assume moreover that the entries of $\partial_h$ come from $E_1$ while the entries of $\partial_v$ come from $E_2$. Does it follow that $\sF$ is a direct sum of box products?
%More precisely, suppose the matrices which define the differential of the Tate resolution $\bT(\sF)$ have no mixed  terms, i.e. all entries 
%are contained in $E_1\tensor_K K  \cup K\tensor_K E_2 \subset E$.
%Does it follow, 
One can ask a similar question for any $t\geq 2$.
\end{question}

\begin{question}
What is the geometric meaning of the other exact subquotient complexes, say the ``half plane'' complexes $\Tate {{c}}  {\{k\}} \o \o$ or $\Tate {{c}}  \o \o {\{k\}}$, defined in Section~\ref{exactness properties}?
\end{question}

\begin{question}
We showed in ~\cite{ES2010} that every complex on $\AA^{m}$
is the direct image of a vector bundle on $\AA^{m}\times \P^{n}$ for some $n$.
Does every complex on $\P^{m}$ occur as a pushforward of a vector bundle
on $\P^{m}\times \PPn$ for some $(n_1,\dots,n_t)$?
\end{question}

\bigskip

\vbox{\noindent Author Addresses:\par
\smallskip
\noindent{David Eisenbud}\par
\noindent{Department of Mathematics, University of California, Berkeley,
Berkeley CA 94720}\par
\noindent{eisenbud@math.berkeley.edu}\par
\smallskip
\noindent{Daniel Erman}\par
\noindent{Department of Mathematics, University of Wisconsin,
Madison WI 53706}\par
\noindent{derman@math.wisc.edu}\par
\smallskip
\noindent{Frank-Olaf Schreyer}\par
\noindent{Mathematik und Informatik, Universit\"at des Saarlandes, Campus E2 4, 
D-66123 Saarbr\"ucken, Germany}\par
\noindent{schreyer@math.uni-sb.de}\par
}

\end{document}